\documentclass[12pt]{amsart}
\usepackage{amsmath}
\usepackage{mathtools}

\DeclarePairedDelimiter{\norm}{\lVert}{\rVert}

\usepackage{amsfonts}
\usepackage{tabularx}
\usepackage{amsthm}
\usepackage{autobreak}
\usepackage[english]{babel}
\usepackage{bold-extra}
\usepackage{hyperref}
\usepackage{makecell}
\usepackage{tikz}
\newtheorem{defi}{Definition}
\newtheorem*{notations}{Notation}

\newtheorem{prop}[defi]{Proposition}
\newtheorem{cor}[defi]{Corollary}

\newtheorem{lem}[defi]{Lemma}

\newtheorem{thm}[defi]{Theorem}

\newtheorem{examples}[defi]{Examples}

\newtheorem{rem}[defi]{Remark}
\numberwithin{defi}{section}
\newcommand{\thistheoremname}{}
\newtheorem*{genericthm*}{\thistheoremname}
\newenvironment{namedthm*}[1]
{\renewcommand{\thistheoremname}{#1}%
	\begin{genericthm*}}
	{\end{genericthm*}}
\allowdisplaybreaks

\def\D{{\mathbb{D}}}

\def\R{{\mathbb{R}}}
\def\C{{\mathbb{C}}}

\def\N{{\mathbb{N}}}
\def\Z{{\mathbb{Z}}}

\def\DD{{\mathcal{D}}}

\def\PP{{\mathcal{P}}}

\def\AA{{\mathcal{A}}}

\newcommand{\la}{\ell_a}

\usepackage{multicol}
\usepackage[T1]{fontenc}
\usepackage[utf8]{inputenc}
\usepackage[a4paper]{geometry}
\usepackage{amssymb,amsmath}
\usepackage{version}
\usepackage{mathtools} 
\usepackage{enumerate}
\usepackage{dsfont}
\usepackage{mathdots}
\makeatletter

\usepackage{latexsym}
\usepackage{color}
\usepackage[normalem]{ulem}

\newenvironment{blueadd}{\begingroup}{\endgroup}

\definecolor{skyblue}{rgb}{0.0,0.55,0.85}

\definecolor{darkgreen}{rgb}{0,0.38,0}

\definecolor{mybrown}{rgb}{0.59,0.29,0.0}

\begin{document}

\title[Composition operators on Beurling--Sobolev spaces]{Composition operators on weighted Beurling--Sobolev spaces}

\author{St\'{e}phane Charpentier}
\address{Aix-Marseille Univ, CNRS, I2M, Marseille, France}
\email{stephane.charpentier.1@univ-amu.fr}

\author{Rachid Zarouf}
\address{Aix Marseille Univ, Lab ADEF, Campus Univ St Jerome, 52 Ave Escadrille Normandie, F-13013 Marseille, France}
\address{Aix-Marseille Univ., Univ. de Toulon, CNRS, CPT, Marseille, France}

\email{rachid.zarouf@univ-amu.fr}

\thanks{The first author is supported by the grant ANR-24-CE40-0892-01 of the French National Research Agency
	ANR (project COMOP). The work of Rachid Zarouf is supported by the pilot center Ampiric, funded by the France 2030 Investment Program operated by the Caisse des D\'ep\^ots.}


\begin{abstract}
Given a function $\varphi:\D\to \D$ holomorphic in a neighborhood of the closed unit disk $\overline{\D}$, we study the matrix
\[
A_\varphi = \big(\widehat{\varphi^n}(k)\big)_{n,k\geq 0}.
\]
We obtain a complete characterization of those $\varphi$ for which $A_\varphi$ acts boundedly from $\ell^p$ to $\ell^q$, $p,q\in [1,\infty]$. These results can be reformulated in terms of composition operators, providing a complete description of the symbols $\varphi$ holomorphic in a neighbourhood of $\overline{\D}$ that induce a bounded composition operator $C_\varphi:f\mapsto f\circ \varphi$ between two possibly different Beurling--Sobolev spaces.

In particular, we prove that $C_\varphi$ is bounded on $\ell_a^p$ for some $p\neq 2$ if and only if it is bounded on every $\ell_a^p$, $p\in [1,\infty]$, if and only if $\varphi$ satisfies an analytic condition at the points of $\overline{\D}$ where it touches the unit circle.

Extensions to weighted Beurling--Sobolev spaces are also provided.
\end{abstract}

\maketitle

\section{Introduction, main results and outline of the proofs}

Let $\D:=\{z\in \C:\, |z|<1\}$ denote the unit disk of the complex plane and let $H(\D)$ be the Fr\'echet space of analytic functions on $\D$. We denote by $H(\overline{\D})$ the subspace of $H(\D)$ consisting of all functions $f$ which extend holomorphically beyond $\overline{\D}$. A function $\varphi \in H(\overline{\D})$ is called a symbol if it is a nonconstant self-map of $\D$.

\medskip

For a symbol $\varphi \in H(\overline{\D})$ and $n,k\in \N$, let
\[
\widehat{\varphi^n}(k)=\frac{1}{2\pi}\int_0^{2\pi}\varphi^n(e^{i t})e^{-ikt}dt
\]
be the $k$th Fourier coefficient of $\varphi^n$.
The central object of this paper is the matrix
\begin{equation*}\label{matrixAphi}
	A_\varphi := \left( \widehat{\varphi^n}(k) \right)_{n,k\geq 0}.
\end{equation*}
We are interested in the problem of describing those $p\in [1,\infty]$ and $\varphi\in H(\overline{\D})$ for which
\[
A_\varphi : \ell^p(w) \to \ell^p(w)
\]
is bounded. Here and below, $w=(w_k)_{k\geq 0}$ is a weight, that is a sequence of positive real numbers satisfying
\begin{equation*}\label{condition-0-weight}
	\liminf_{k\to \infty}w_k^{1/k}\geq 1,
\end{equation*}
and, for $p\in [1,\infty)$,
\[
\ell^p(w):=\{a=(a_n)_{n\geq 0}:\, \Vert a \Vert_p^p:=\sum_{k\geq 0}|a_k|^pw_k^p < \infty\}.
\]
We also denote
\[
\ell^\infty(w):=\{a=(a_n)_{n\geq 0}:\, \Vert a \Vert_\infty:=\sup_{k\geq 0}|a_k|w_k < \infty\}.
\]

This problem is closely related to that of describing bounded composition operators acting on the analytic versions of the $\ell^p(w)$ spaces. Given a symbol $\varphi$ and a space $X$ of analytic functions on $\D$, the composition operator of symbol $\varphi$ is defined as
\[
C_\varphi: f\mapsto f\circ \varphi,\quad f\in X.
\]
Then, for $p\in[1,\infty]$ and  a weight $w$, let
\[
\ell_{a}^{p}(w):= \{f=\sum_{k\geq 0}a_kz^k\in H(\D):\, (a_k)_{k\geq 0} \in \ell^p(w)\}.
\]
The infinite matrix of $C_\varphi$, written in the canonical basis $(z^n)_{n\geq0}$ of $\ell_a^p(w)$, coincides with the transpose $A_\varphi^t$ of $A_\varphi$. For example, if $w=1$, then the boundedness of $C_\varphi$ on $\ell_a^p(w)$ is equivalent to the boundedness of the matrix $A_\varphi^t$ as an operator acting on $\ell^p(w)$.

\medskip

The spaces $\ell_a^p(w)$ are often referred to as the weighted Beurling--Sobolev spaces. For every weight $w$ and every $p\in[1,\infty]$, the space $\ell_a^p(w)$, $p\in [1,\infty]$, is continuously embedded into $H(\D)$. Two weights $w$ and $\beta$ are said to be comparable, and we write $w\asymp \beta$, if there exist constants $c,C>0$ such that
\[
c \beta _n \leq w_n \leq C \beta_n,\quad n\in \N.
\]
In this case, $\ell_a^p(w)=\ell_a^p(\beta)$ with equivalent norms. When $w$ is comparable to a constant, we simply write $\ell_a^p$ instead of $\ell_a^p(w)$ (and $\ell^p$ instead of $\ell^p(w)$). Classical examples of Beurling--Sobolev spaces are the Wiener space $W:=\ell_a^1$, the Hardy space $H^2:=\ell_a^2$, the Bergman space $A^2:=\ell_a^2(w)$ for $w\asymp (n+1)^{-1/2}$, and the Dirichlet space $\DD:=\ell_a^2(w)$ for $w\asymp (n+1)^{1/2}$.

Composition operators have been extensively studied on the Hardy, Bergman and Dirichlet spaces and, more generally, on Banach spaces of analytic functions whose norm can be represented as an integral, allowing changes of variables. For example, it is well-known that every composition operator is bounded on $H^2$ and $A^2$, while some symbols induce an unbounded composition operator on the Dirichlet space. See \cite{ShapiroBook1993} for a complete reference on this classical subject. For general weights $w$, or when $p\neq 2$, such integral representations are no longer available in $\ell_a^p(w)$, making the problem substantially more difficult and, comparatively, little is known about composition operators on these spaces.

In her thesis, N. Zorboska explicitly stated the problem of characterizing those weights $w$ for which every symbol $\varphi \in \ell_a^2(w)$ induces a composition operator mapping $\ell_a^2(w)$ into itself \cite{Zorboska1988}. In 2025, Lef\`evre, Li, Quéffelec and Rodr\'iguez-Piazza solved this problem, by proving that every composition operator is bounded on $\ell_a^2(w)$ if and only if $w$ is comparable to a decreasing weight and is slowly oscillating \cite{LLQR-l2weight-2025}. We say that a weight $w$ is slowly oscillating if, for every $0<J<1<L<\infty$,
\[
w_k\asymp w_n, \quad Jn \leq k\leq Ln.
\]
In passing, they established that the composition operators $C_\varphi$ are bounded on $\ell_a^2(w)$ for all $\varphi \in H(\overline{\D})$ if and only if $w$ is slowly oscillating. In \cite{LLQR-lpa-2025}, the same authors proved that this phenomenon is specific to $p=2$: if $p \neq 2$, no nontrivial automorphism induces a bounded composition operator on $\ell_a^p(w)$, for any weight $w$.

A different line of results goes back to 1969 \cite{Newman1969}. In this paper, Newman characterized those symbols $\varphi$ in $H(\overline{\D})$, for which $C_\varphi$ is bounded on the Wiener algebra $\ell_a^1$. In 1993, Blyudze and Shimorin \cite{BlyudzeShimorin1993} further characterized those $\varphi \in H(\overline{\D})$ for which $C_\varphi$ acts boundedly from $\ell_a^1$ to $\ell_a^p$. It is not difficult to see that the boundedness of $C_\varphi$ from $\ell_a^1$ to $\ell_a^p$ is equivalent to a uniform control on the $\ell^p$ norm of the rows of the matrix $A_\varphi$: $C_\varphi$ is bounded from $\ell_a^1$ to $\ell_a^p$ if and only if
\begin{equation}\label{CB-boundedness-nec-cond}
	\sup_{n\geq 0}\Vert \varphi^n \Vert_{\ell_a^p}^p=\sup_{n\geq 0} \sum_{k\geq 0}|\widehat{\varphi^n}(k)|^p <\infty.
\end{equation}
Since $A_\varphi =C_\varphi^t$, we note in passing that \eqref{CB-boundedness-nec-cond} is also equivalent to $A_\varphi:\ell^{p'} \to \ell^\infty$ bounded, where $p'$ is the conjugate exponent of $p$ (\textit{i.e.} $1/p+1/p'=1$). Blyudze and Shimorin established sharp estimates of $\Vert \varphi ^n \Vert _{\ell_a^p}$ (see Theorem~\ref{thm-asympt-BS} below), allowing them to describe those symbols $\varphi \in H(\overline{\D})$ for which \eqref{CB-boundedness-nec-cond} holds.

\medskip

In this paper, we follow this second line of research. In order to state our results, we need to introduce some notation from \cite{BlyudzeShimorin1993,Newman1969}, that will be used repeatedly in the sequel. For a symbol $\varphi$ in $H(\overline{\D})$, set
\[
J_{\varphi}:=\{\theta\in[0,2\pi): |\varphi(e^{i\theta})|=1\}.
\]
One can easily see that if $J_{\varphi} =\emptyset$, then $\varphi(\overline{\D}) \subset \D$ and then $C_\varphi$ (and hence $A_\varphi$ as well) is bounded on $\ell_a^p$ for every $p\in [1,\infty]$. When $J_\varphi\neq \emptyset$, and for $\theta\in J_{\varphi}$, let us write the expansion
\[
\log \varphi(e^{i(\theta+t)})=\sum_{j\ge 0} c_j(\theta)t^j
\]
for $t$ small enough, and define
\[
r(\theta):=\inf\{j\ge 2:\ c_j(\theta)\neq 0\}
\quad \text{and} \quad
\sigma(\theta):=\inf\{j\ge 2:\ \Re c_j(\theta)\neq 0\}.
\]
Here and below, we use the convention $\inf\emptyset = \infty$. Clearly, $r(\theta) \leq \sigma(\theta)$ for every $\varphi$ and every $\theta \in J_\varphi$. Moreover, one can check that $\sigma(\theta) = \infty$ for some $\theta \in J_\varphi$ if and only if $\varphi$ is a finite Blaschke product, while $r(\theta)=\infty$ for some $\theta \in J_\varphi$ if and only if $\varphi$ is monomial; see \cite{Hedstrom1966} or \cite[Lemma 1.1]{BlyudzeShimorin1993}.

For $p\in[1,\infty)$, set
\[
	\alpha_p(\varphi):=
	\max_{\theta\in J_{\varphi}}
	\frac{(\sigma(\theta)-r(\theta))(2-p)-(2p-2)}{2p\,\sigma(\theta)}
\]
if $\varphi$ is not a finite Blaschke product, and
\[
	\alpha_p(\varphi):= \frac{2-p}{2p}
\]
if $\varphi$ is a non-monomial finite Blaschke product. A symbol $\varphi$ is called a monomial if $\varphi =\mu z^n$ for some $|\mu|\leq 1$ and $n\in \N\setminus \{0\}$. When $\varphi$ is monomial, the problem we are interested in becomes trivial.

The result of Blyudze and Shimorin can then be stated as follows.

\begin{thm}[Corollary 1.10 in \cite{BlyudzeShimorin1993}]\label{BZ-1993}Let $p\in [1,\infty)$ and let $\varphi \in H(\overline{\D})$ non-monomial such that $J_\varphi \neq \emptyset$. The following assertions are equivalent:
	\begin{enumerate}
		\item $\sup_{n\geq 0}\Vert \varphi^n \Vert_{\ell_a^p} < \infty$;
		\item $A_\varphi:\,\ell^{p'} \to \ell^\infty$ is bounded;
		\item $C_\varphi:\,\ell_a^1 \to \ell_a^p$ is bounded;
		\item $\alpha_p(\varphi)\leq 0$.
	\end{enumerate}
\end{thm}

We point out that the problem of estimating the $\ell^p$ norms, $p\in[1,\infty]$, of the powers of a function $\varphi \in H(\D)$ is classical and has numerous applications in, for example, harmonic analysis, operator theory, complex function theory, or the study of crystalline measures; see \cite{BFZ2024,BFZ-Blaschkeproduct-2024,Kahane1956,LLQR-l2weight-2025,Meyer2021,SzehrZarouf2026Crystalline,SzehrZarouf2020,SzehrZarouf2021}.
\medskip

Our first aim is to obtain a statement similar to Theorem~\ref{BZ-1993} for the columns of $A_\varphi$, namely, a characterization of those $\varphi \in H(\overline{\D})$ for which
\begin{equation*}\label{column-CNS}
	\sup_{k\geq 0}\Vert (\widehat{\varphi^n}(k))_{n\geq 0} \Vert_{\ell^p}^p =\sup_{k\geq 0} \sum_{n\geq 0}|\widehat{\varphi^n}(k)|^p <\infty.
\end{equation*}

After observing that
\[
C_\varphi: \ell_a^{p'}\to \ell_a^\infty \text{ bounded}  \quad \Longleftrightarrow \quad  \sup_{k\geq 0} \sum_{n\geq 0}|\widehat{\varphi^n}(k)|^p <\infty \quad \Longleftrightarrow \quad A_\varphi: \ell^1 \to \ell^p\text{ bounded},
\]
we establish the following.
\begin{thm}\label{thm-intro-main}Let $p\in [1,\infty)$ and let $\varphi \in H(\overline{\D})$ non-monomial such that $J_\varphi \neq \emptyset$. The following assertions are equivalent:
	\begin{enumerate}
		\item $\sup_{k\geq 0}\Vert (\widehat{\varphi^n}(k))_{n\geq 0} \Vert_{\ell^p} < \infty$;
		\item $A_\varphi:\,\ell^1 \to \ell^p$ is bounded;
		\item $C_\varphi:\,\ell_a^{p'} \to \ell_a^\infty$ is bounded;
		\item $\alpha_p(\varphi)\leq 0$.
	\end{enumerate}
\end{thm}

Theorem~\ref{thm-intro-main} has several consequences, some of which are rather unexpected, especially when formulated in terms of composition operators (the same consequences are valid for $A_\varphi$). The first one is a direct consequence of Theorems~\ref{BZ-1993} and \ref{thm-intro-main}.

\begin{cor}\label{main-thm-intro-1}Let $p\in [1,\infty)$ and let $\varphi \in H(\overline{\D})$ non-monomial such that $J_\varphi \neq \emptyset$. The following assertions are equivalent:
	\begin{enumerate}
		\item $\sup_{k\geq 0}\Vert (\widehat{\varphi^n}(k))_{n\geq 0} \Vert_{\ell^p}  < \infty$;
		\item $\sup_{n\geq 0}\Vert (\widehat{\varphi^n}(k))_{k\geq 0} \Vert_{\ell^p} < \infty$;
		\item $C_\varphi$ maps $\ell_a^1$ into $\ell_a^p$;
		\item $C_\varphi$ maps $\ell_a^{p'}$ into $\ell_a^\infty$;
		\item $\alpha_p(\varphi)\leq 0$.
	\end{enumerate}
\end{cor}
The equivalence between (1) and (2) is a priori highly non-trivial and reflects the special structure of $A_\varphi$. Since $\alpha_p(B) \leq 0$ if and only if $p\geq 2$ when $B$ is a finite Blaschke product, we deduce:

\begin{cor}\label{cor-BF-intro}For every non-monomial finite Blaschke product $B$,
	\[
	C_B:\, \ell_a^1\to \ell_a^p\text{ bounded} \quad \Longleftrightarrow \quad p\geq 2,
	\]
	and
	\[
	C_B:\, \ell_a^p\to \ell_a^\infty\text{ bounded} \quad \Longleftrightarrow  \quad p\leq 2.
	\]
	In particular, $C_B$ is bounded on $\ell_a^p$ if and only if $p= 2$.
\end{cor}

Taking $p=1$ in Corollary~\ref{main-thm-intro-1}, we deduce:

\begin{cor}\label{cor-linfty-intro}Let $\varphi \in H(\overline{\D})$ non-monomial such that $J_\varphi \neq \emptyset$. The following assertions are equivalent:
	\begin{enumerate}
		\item $C_\varphi$ maps $\ell_a^{1}$ into itself;
		\item $C_\varphi$ maps $\ell_a^{\infty}$ into itself;
		\item $C_\varphi$ maps $\ell_a^{p}$ into itself for every $p\in[1,\infty]$;
		\item $r(\theta)=\sigma(\theta) <\infty $ for every $\theta \in J_\varphi$.
	\end{enumerate}
\end{cor}

\medskip

The previous results are solutions to the initial problem of the boundedness of $C_\varphi$ acting on $\ell_a^p$, essentially for $p=1$ and $p=\infty$. Indeed, a uniform control on the $\ell^p$ norms of the rows and the columns of a matrix need not imply the boundedness of the matrix on $\ell^p$, $p\in[1,\infty]$; except if such a control holds for the extremities $p=1$ or $p=\infty$, as in Corollary~\ref{cor-linfty-intro}, in which can one strikingly get a control on the operator norm of $C_\varphi$ for every $p$. To provide a complete solution to our initial problem, it is thus necessary to understand the behaviour of $C_\varphi$ on $\ell^p_a$ (or $A_\varphi$ on $\ell^p$) when $r(\theta)<\sigma(\theta)$ for some $\theta \in J_\varphi$, for \textit{any} $p\in [1,\infty]$. This is achieved by our second main result:

\begin{thm}\label{main-thm-rlesssigma}Let $p\in [1,\infty]$ and let $\varphi \in H(\overline{\D})$ non-monomial such that $J_\varphi \neq \emptyset$. If $r(\theta)<\sigma(\theta)<\infty$ for some $\theta\in J_\varphi$, then $C_\varphi: \ell_a^p \to \ell_a^p$ is bounded if and only if $p=2$.
\end{thm}

Altogether, Corollary~\ref{cor-linfty-intro} and Theorem~\ref{main-thm-rlesssigma} therefore provide a complete description of those $\varphi\in H(\overline{\D})$ for which $C_\varphi$ is bounded on $\ell_a^p$, for every fixed $p\in [1,\infty]$, solving the initially stated problem. 

Furthermore, these two results and the proof of the latter theorem turn out to be sufficiently powerful to allow us, with the help of the Riesz--Thorin interpolation theorem, to even establish a full description of those symbols $\varphi \in H(\overline{\D})$ and those exponents $p,q\in  [1,\infty]$ for which $C_\varphi: \ell_a^p \to \ell_a^q$ is bounded. This description can be stated as follows.

\begin{thm}\label{Thm-Fulldescription}Let $\varphi\in H(\overline{\D})$ be a symbol. Then the following assertions hold.
	\begin{enumerate}
		\item If $J_\varphi=\varnothing$, then $C_\varphi:\ell_a^p\to\ell_a^q$ is bounded for every $p,q\in[1,\infty]$.
		\item If $\varphi = \mu z^n$ is a monomial, for some $|\mu|\leq 1$ and some $n\geq 1$, then
		\begin{enumerate}
			\item If $|\mu|< 1$, then $C_\varphi:\ell_a^p\to\ell_a^q$ is bounded for every $p,q\in[1,\infty]$;
			\item If $|\mu|=1$, then, for every $p,q\in[1,\infty]$, $C_\varphi:\ell_a^p\to\ell_a^q$ is bounded if and only if $p\leq q$.
		\end{enumerate}
		\item If $\varphi$ is a non-monomial finite Blaschke product, then $C_\varphi:\ell_a^p\to\ell_a^q$ is bounded if and only if $1\leq p \leq 2 \leq q \leq \infty$.
		\item Assume that $J_\varphi\neq\varnothing$ and that $\varphi$ is neither a monomial nor a finite Blaschke product. Set
		\[
		D_\varphi:=\max_{\theta\in J_\varphi}	\left(\sigma(\theta)-r(\theta)\right).
		\]
		Then, for every $p,q\in[1,\infty]$, the following assertions are equivalent:
		\begin{enumerate}
			\item $C_\varphi:\ell_a^p\to\ell_a^q$ is bounded;
			\item The conditions
			\[
			(D_\varphi+1)\frac1p-\frac1q -\frac{D_\varphi}{2}\geq 0 \quad \text{and}\quad \frac1p-(D_\varphi+1)\frac1q+ \frac{D_\varphi}{2}\geq 0
			\]
			hold, with the convention $1/\infty =0$.
		\end{enumerate}
		In particular, if $D_\varphi=0$, then (b) is equivalent to $p\leq q$.
	\end{enumerate}
\end{thm}

Let us say a word on the proofs of Corollary~\ref{cor-linfty-intro} and Theorem~\ref{main-thm-rlesssigma}. Both of them rely on some precise estimates of the Fourier coefficients of $\varphi^n$ but the purposes are different. For Corollary~\ref{cor-linfty-intro}, one has to upper bound the $\ell^p$ norms of the columns of the matrix $A_\varphi^t$ while, for Theorem~\ref{main-thm-rlesssigma}, one has to extract from $A_\varphi^t$ a sequence of finite dimensional block matrices that is unbounded on finite dimensional $\ell^p$ spaces. For the first purpose, we need to upper bound the integral representing the Fourier coefficient $\widehat{\varphi^n}(k)$. In the case where $\varphi$ is not a finite Blaschke product, which is the most delicate one, the main contributions in this integral are due to the finitely many points where $\varphi$ touches the boundary of $\D$, that is near the points of $J_\varphi$. The upper estimate of the Fourier integral localized near the points $\theta \in J_\varphi$ will be derived from a van der Corput's type lemma, whose intricate verification of the assumptions requires to distinguish the case $r(\theta)=\sigma(\theta)$ and $r(\theta) < \sigma(\theta)$. The proof of Theorem~\ref{main-thm-rlesssigma} consists in obtaining sharp estimates for the contribution of sufficiently many Fourier coefficients $\widehat{\varphi^n}(k)$ near a fixed point of $J_\varphi$, in order to define finite dimensional subblocks of $A_\varphi^t$ whose operator norms grow unboundedly. The crucial ingredient is a uniform version of the stationary phase method, due to Fedoryuk \cite{Fedoryuk1977}.

\medskip

We point out that the quantities $r(\theta)$ and $\sigma(\theta)$ are rather handy. In particular, it is not difficult to provide examples of symbols $\varphi \in H(\overline{\D})$ that satisfy $D_\varphi >0$ or $D_\varphi=0$; see Examples~\ref{examplesymbols}.

\medskip

In a second part of the paper, we shall investigate the validity of the previous results in the weighted setting. For simplicity, we only consider the same exponents $p=q\in [1,\infty]$. We establish that all the previous results extend, in this setting, whenever the weight is slowly oscillating, namely:
\begin{thm} Let $w$ be a slowly oscillating weight and let $p\in[1,\infty]$. In Corollaries~\ref{main-thm-intro-1}--\ref{cor-linfty-intro} and Theorems~\ref{main-thm-rlesssigma} and \ref{Thm-Fulldescription} (with $p=q$), one can replace $\ell_a^p$ (resp. $\ell^p$) with $\ell_a^p(w)$ (resp. $\ell^p(w)$).
\end{thm}

We will give counter-examples showing that, when the slow oscillation is dropped, none of the above results hold in general.

\medskip

In all the previous, $C_\varphi$ can be replaced with $A_\varphi$ and $\ell_a^p$ by $\ell^p$. The matrix $A_\varphi$ is a rather natural object that can appear in many problems. Thus, the results of this paper and the techniques used in the corresponding proofs may have a broader interest in Analysis.

\medskip

The paper is organized as follows. Section~\ref{sec-prelim} contains the preliminary results. Section~\ref{sec-row-column} is devoted to the proof of Theorem~\ref{thm-intro-main}, while Theorems~\ref{main-thm-rlesssigma} and \ref{Thm-Fulldescription} are proved in Section~\ref{sec-thmfull}. The weighted setting is discussed in the final section.

\begin{notations}{\rm Throughout the paper, we use the standard Landau notation $O(\cdot)$ and $o(\cdot)$. For nonnegative quantities $A$ and $B$, we write $A\lesssim B$ if there exists a constant $C>0$ such that $A\leq CB$, and $A\asymp B$ if both $A\lesssim B$ and $B\lesssim A$. The implicit constants may depend on fixed parameters.
	
	We shall consider symbols $\varphi \in H(\overline{\D})$. An important subclass of $H(\overline{\D})$ consists of finite Blaschke products. We recall that, given a finite sequence $\Lambda=(\lambda_1,\ldots,\lambda_d)\subset \D$ and a unimodular constant $\mu$, the finite Blaschke product $B$ associated with $\Lambda$ is defined by
	\begin{equation*}\label{Blaschke-def}
		B(z)=\mu\prod_{j=1}^d\frac{\lambda_j -z}{1-\overline{\lambda_j}z},\quad z\in \D.
	\end{equation*}
	The natural number $d$ is called the degree of $B$. When $d=1$, $B$ is called a single Blaschke factor. We shall use the notation $b_{\lambda}$ for the single Blaschke factor vanishing at $\lambda$. Recall that the automorphisms of $\D$ are exactly the single Blaschke factors. We shall also call a monomial any function of the form $\mu z^n$ for some $\mu \in \overline{\D}$ and some $n\in \N\setminus\{0\}$. }
\end{notations}

\section{Preliminary results}\label{sec-prelim}

This section contains preliminary or elementary useful results. The first lemma will be used repeatedly in several forthcoming proofs.

\begin{lem}\label{lem-elementary}Let $p\in [1,\infty)$ and let $\varphi \in H(\overline{\D})$ be  a symbol. Then there exist constants $0<J<1<L<\infty$ such that, for every $\gamma \geq 0$, there exists $s\in(0,1)$ such that, for every $n$ large enough,
	\[
	\sum_{k<Jn}\left|\widehat{\varphi^{n}}(k)\right|^{p}k^\gamma+\sum_{k>Ln}\left|\widehat{\varphi^{n}}(k)\right|^{p}k^\gamma \leq s^n
	\]
	and, for every $k$ large enough,
	\[
	\sum_{n<Jk}\left|\widehat{\varphi^{n}}(k)\right|^{p} n^\gamma + \sum_{n>Lk}\left|\widehat{\varphi^{n}}(k)\right|^{p}n^\gamma \leq s^k.
	\]
\end{lem}

\begin{proof}Let $\varrho>1$ be such that $\varphi$ is analytic in the open disk $D(0,\varrho)$. Fix $0<r<1$, $1<R<\varrho$ and set
	\[
	M_r:=\max_{|z|=r}|\varphi(z)|<1\quad \text{and} \quad M_R:=1+\max_{|z|=R}|\varphi(z)|>1.
	\]
	Since $M_r < 1$ and $M_R >1$, there exist two constants $0<J<1<L$ such that
	\begin{equation}\label{eq-lem-elem}
		M_r<r^J,\,M_R<R^L\quad  \text{and}\quad M_R^{J}<R,\, M_r^{L}<r.
	\end{equation}
	Applying Cauchy's estimates to $\varphi^{n}$ on the circles $|z|=r$ and $|z|=R$, and then summing over $k\in [0,Jn)$ and over $k>Ln$, respectively, we obtain, for every $n\ge 0$,
	\[
	\sum_{k< Jn}\left|\widehat{\varphi^{n}}(k)\right|^{p} \leq \frac{1}{1-r^p}\left(\frac{M_r}{r^J}\right)^{pn} \quad \text{and} \quad \sum_{k> Ln}\left|\widehat{\varphi^{n}}(k)\right|^{p} \leq \frac{R^p}{R^p-1}\left(\frac{M_R}{R^L}\right)^{pn}.
	\]
	Similarly, applying now Cauchy's estimates to $\varphi^{n}$ on the circles $|z|=R$ and $|z|=r$, and summing over $n\in [0,Jk]$ and over $n\geq Lk$, respectively, we obtain, for every $k\ge 0$,
	\[
	\sum_{n<Jk}\left|\widehat{\varphi^{n}}(k)\right|^{p} \leq \frac{M_R^p}{M_R^p-1}\left(\frac{M_R^J}{R}\right)^{pk} \quad \text{and} \quad \sum_{n>Lk}\left|\widehat{\varphi^{n}}(k)\right|^{p} \leq \frac{1}{1-M_r^p}\left(\frac{M_r^L}{r}\right)^{pk}.
	\]
	By  \eqref{eq-lem-elem},
	\[
	\max\left(\frac{M_r}{r^J},\frac{M_R}{R^L},\frac{M_R^J}{R},\frac{M_r^L}{r}\right)<1,
	\]
	which guarantees the existence of the desired $s\in (0,1)$ for $\gamma=0$; the proof of the case $\gamma>0$ is similar.
\end{proof}

Note that, as in the proof of Lemma~\ref{lem-elementary}, Cauchy's estimates also yield the following. The details are left to the reader.

\begin{prop}\label{propCauchyesti}Let $\varphi \in H(\overline{\D})$ be a symbol. There exists $s\in (0,1)$ such that,
	\begin{enumerate}
		\item for every $n\geq 0$ and every $k$ large enough (depending on $n$), one has $|\widehat{\varphi^n}(k)|\leq s^k$;
		\item for every $k\geq 0$ and every $n$ large enough (depending on $k$), one has $|\widehat{\varphi^n}(k)|\leq s^n$.
	\end{enumerate}
	In particular, for every $p\in [1,\infty]$ and every weight $w$,
	\[
	\Vert (\widehat{\varphi^n}(k))_{k\geq0}\Vert_{\ell^p(w^{-1})} < \infty,\,\,n\geq 0 \quad \text{and}\quad \Vert(\widehat{\varphi^n}(k))_{n\geq 0}\Vert_{\ell^p(w^{-1})}< \infty,\,\,k\geq 0.
	\]
\end{prop}

The next result contains the equivalences (1)$\Longleftrightarrow$(2)$\Longleftrightarrow$(3) in both Theorems~\ref{BZ-1993} and \ref{thm-intro-main}.

\begin{prop}\label{prop-key-boundedness-carac}Let $w$ be a weight and let $\varphi \in H(\overline{\D})$ be a symbol. Then
	\begin{enumerate}\item\label{item1prop} For every $p\in[1,\infty)$,
		\[
		\begin{array}{lll}
			A_\varphi: \ell^{p'}(w^{-1}) \to \ell^\infty(w^{-1})\text{ bounded} & \Longleftrightarrow & C_\varphi: \ell_a^{1}(w)\to \ell_a^p(w) \text{ bounded}\\
			& \Longleftrightarrow & \displaystyle{\sup_{n\geq 0} \sum_{k\geq 0}|\widehat{\varphi^n}(k)|^p\frac{w_k^p}{w_n^p} <\infty }.
		\end{array}
		\]
		\item\label{item2prop} For every $p\in(1,\infty]$,
		\[
		\begin{array}{lll}
			A_\varphi: \ell^1(w^{-1}) \to \ell^{p'}(w^{-1})\text{ bounded} & \Longleftrightarrow &  C_\varphi: \ell_a^{p}(w)\to \ell_a^\infty(w) \text{ bounded}\\
			& \Longleftrightarrow & \displaystyle{\sup_{k\geq 0} \sum_{n\geq 0} |\widehat{\varphi^n}(k)|^{p'}\frac{w_k^{p'}}{w_n^{p'}} <\infty }.
		\end{array}
		\]
	\end{enumerate}
\end{prop}

\begin{proof}We only prove the second implication in both \eqref{item1prop} and \eqref{item2prop}. Since $A_\varphi$ is the transpose of the matrix of $C_\varphi$ written in the canonical basis of $\ell_a^p(w)$, the proofs of the two other implications are very similar.
				
For $f(z)=\sum_{n\geq 0}a_nz^n \in \ell_a^\infty(w)$, we have
	\begin{equation}\label{eq-prop1-proof}
		f\circ \varphi(z)=\sum_{n\geq 0}a_n\sum_{k\geq 0} \widehat{\varphi^n}(k)z^k=\sum_{k\geq 0}\left(\sum_{n\geq 0} a_n\widehat{\varphi^n}(k)\right)z^k,\quad z\in \D.
	\end{equation}
	Notice that, by Proposition~\ref{propCauchyesti}, all the involved series are convergent.
	
	For the second equivalence in Assertion~(1), remark first that
	\begin{eqnarray*}
	\Vert f\circ \varphi \Vert _{\ell_a^p(w)} & \leq & \sum_{n\geq 0}|a_n|w_n \Vert \varphi^n \Vert _{\ell_a^p(w)} w_n^{-1},\\
	& \leq & \Vert f \Vert _{\ell_a^1(w)}\sup_{n\geq 0}\left(\sum_{k\geq 0} |\widehat{\varphi^n}(k)|^p\frac{w_k^p}{w_n^p}\right)^{1/p},
	\end{eqnarray*}
 	hence
 	\[
 	\sup_{n\geq 0}\sum_{k\geq 0} |\widehat{\varphi^n}(k)|^p\frac{w_k^p}{w_n^p} < \infty \, \Longrightarrow \, C_\varphi:\ell_a^1(w) \to \ell_a^p(w)\text{ bounded}.
 	\]
	
	The converse follows from the observation
	\[
	\sum_{k\geq 0} |\widehat{\varphi^n}(k)|^p\frac{w_k^p}{w_n^p} = \frac{\Vert \varphi ^n \Vert^p _{\ell_a^p(w)}}{\Vert z^n \Vert _{\ell_a^1(w)}^p}.
	\]
	
	\medskip
	
	Let us prove the second equivalence in Assertion~\eqref{item2prop}. Let us first assume that $p < \infty$. For every fixed $k\geq 0$, consider the linear functional $T_k:\ell_a^p(w) \to \C$ defined by
	\[
	T_k(f)=\sum_{n\geq 0}a_n\widehat{\varphi^n}(k)w_k,\quad f(z)=\sum_{n\geq 0}a_nz^n.
	\]
	Using \eqref{eq-prop1-proof}, we get
	\[
	C_\varphi:\,\ell_a^p(w) \to \ell_a^\infty(w) \text{ bounded} \quad \Longleftrightarrow \quad \sup_{k\geq 0}\vert T_k(f) \vert<\infty,\, f \in \ell_a^p(w).
	\]
	Moreover, the uniform boundedness principle yields
	\[
	\sup_{k\geq 0}\vert T_k(f) \vert<\infty,\, f \in \ell_a^p(w)  \quad \Longleftrightarrow \quad \sup_{k\geq 0}\Vert T_k\Vert < \infty.
	\]
	Now, recall that the dual space of $\ell_a^p(w)$ identifies with $\ell_a^{p'}(w^{-1})$ via the Cauchy pairing, which agrees with the usual inner product $\left<\cdot,\cdot\right>$ on $\ell_a^2$. Thus, one has
	\[
	\sup_{k\geq 0}\Vert T_k \Vert = \sup_{k\geq 0}\Vert (\widehat{\varphi^n}(k)w_k)_n \Vert _{\ell^{p'}(w^{-1})} = \sup_{k\geq 0}\left(\sum _{n\geq 0} |\widehat{\varphi^n}(k)| ^{p'}\frac{w_k^{p'}}{w_n^{p'}}\right)^{1/p'}.
	\]
	The desired equivalence follows.
	
	To finish, let us assume that $p=\infty$. For every $f(z)=\sum_{n\geq 0}a_nz^n \in \ell_a^\infty(w)$ and every $k\geq 0$, one has
	\[
	|\widehat{f\circ \varphi}(k)|w_k=\left|\sum_{n\geq 0}a_n\widehat{\varphi^n}(k)\right|w_k\leq \sup_{n\geq 0} |a_nw_n|\sum_{n\geq 0}|\widehat{\varphi^n}(k)|\frac{w_k}{w_n}.
	\]
	Thus,
	\[
	\sup_{k\geq 0}\sum_{n\geq 0}|\widehat{\varphi^n}(k)|\frac{w_k}{w_n} < \infty \quad \Longrightarrow \quad C_\varphi:\, \ell_a^\infty(w) \to \ell_a^\infty(w) \text{ bounded}.
	\]
	Conversely, fix $k\in \N$ and $N\in \N$. For $0\leq n \leq N$, choose
	\[
	a_n  = \left\{
	\begin{array}{ll}\frac{\overline{\widehat{\varphi^n}(k)}}{|\widehat{\varphi^n}(k)|}w_n^{-1} & \text{if } \widehat{\varphi^n}(k)\neq 0\\
		w_n^{-1} &  \text{if } \widehat{\varphi^n}(k)= 0
	\end{array}
	\right.
	\]
	and set $f_N(z)=\sum_{n=0}^N a_nz^n$. One has $\Vert f_N \Vert_{\ell_a^\infty(w)}=1$. Therefore, if $C_\varphi: \, \ell_a^\infty(w) \to \ell_a^\infty(w)$ is bounded, then
	\[
	w_k|\widehat{f_N\circ \varphi}(k)| = \sum_{n=0}^N|\widehat{\varphi^n}(k)|\frac{w_k}{w_n}\leq \Vert C_{\varphi} \Vert _{\ell_a^\infty(w)\to \ell_a^\infty(w)}<\infty.
	\]
	Letting $N\to \infty$ and taking the supremum over $k\geq 0$ proves the remaining implication.
\end{proof}

\section{$\ell^p$ norms of the rows and the columns of the matrix $A_\varphi$}\label{sec-row-column}
In this section, we complete the proof of Theorem~\ref{thm-intro-main}. By Proposition~\ref{prop-key-boundedness-carac} applied with $w$ equal to the constant weight $1$, it remains to prove the following, which is the core of the statement.

\begin{thm}\label{thm-main-1-to-prove}For every $p\in [1,\infty)$ and every $\varphi \in H(\overline{\D})$ non-monomial such that $J_\varphi \neq \emptyset$,
	we have
	\begin{equation*}\label{eq-thm1-core}
		\sup_{k\geq 0}\sum _{n\geq 0} |\widehat{\varphi^n}(k)|^{p} <\infty \quad \Longleftrightarrow \quad \alpha_p(\varphi) \leq 0.
	\end{equation*}
\end{thm}

The proof is a consequence of estimates of $\Vert (\widehat{\varphi^n}(k))_{k\geq 0}\Vert _{\ell^p}$ and $\Vert (\widehat{\varphi^n}(k))_{n\geq 0}\Vert _{\ell^p}$, that may be interesting in themselves. For the first quantity, Blyudze and Shimorin established the following:

\begin{thm}[Theorems 1.2 and 1.3 in \cite{BlyudzeShimorin1993}]\label{thm-asympt-BS}Let $\varphi \in H(\overline{\D})$ be a non-monomial symbol such that $J_\varphi \neq \emptyset$. Then,
	\begin{enumerate}
		\item \label{assert-BS-majoration} for every $p\in [1,2]$,
		\[
			\sup_{n\geq 0}(n+1)^{-p\alpha_p(\varphi)}\sum_{k\geq 0} |\widehat{\varphi^n}(k)|^p < \infty ;
		\]
		\item \label{assert-BS-minoration} for every $p\in [1,\infty)$,
		\[
			\inf_{n\geq 0}(n+1)^{-p\alpha_p(\varphi)}\sum_{k\geq 0} |\widehat{\varphi^n}(k)|^p > 0.
		\]
	\end{enumerate}
\end{thm}

\begin{rem}{\rm We indicate that, in Assertion~\ref{assert-BS-majoration}, the restriction $p\in [1,2]$ is not merely technical. Indeed, when $\varphi$ is a single Blaschke factor, the authors of \cite{SzehrZarouf2020} obtained sharp estimates on $\sum_{k\geq 0} |\widehat{\varphi^n}(k)|^p$ for the whole range $p\in[1,\infty]$, showing that, in this case, Assertion~\ref{assert-BS-majoration} holds if and only if $p\in[1,4)$. In particular, Assertion~\ref{assert-BS-minoration} is not sharp in general for $p\geq 4$.}
\end{rem}

This is the following analogue of Theorem~\ref{thm-asympt-BS} for the quantity $\sum_{n\geq 0} |\widehat{\varphi^n}(k)|^p$ that will give us Theorem~\ref{thm-main-1-to-prove}

\begin{thm}\label{thm-asympt-BS-dual-upper}Let $\varphi \in H(\overline{\D})$ be a non-monomial symbol such that $J_\varphi \neq \emptyset$. Then,
	\begin{enumerate}
		\item \label{eq-asympt-BS-sup-majoration}
		for every $p\in [1,2]$,
		\[
			\sup_{k\geq 0}(k+1)^{-p\alpha_p(\varphi)}\sum_{n\geq 0} |\widehat{\varphi^n}(k)|^p < \infty ;
		\]
		\item \label{assertion-minoration-dual} for every $p\in [1,\infty)$, there exist constants $0<J<L<\infty$, $c>0$ and $N_0\in \N$ such that, for every $N\geq N_0$,
		\[
			\frac{1}{N}\sum_{JN\leq k \leq LN}\sum_{n\geq 0} |\widehat{\varphi^n}(k)|^p \geq cN^{p\alpha_p(\varphi)}.
		\]
		In particular, for every $p\in [1,\infty)$,
		\[
			\limsup_{k}(k+1)^{-p\alpha_p(\varphi)}\sum_{n\geq 0} |\widehat{\varphi^n}(k)|^p >0.
		\]
	\end{enumerate}
\end{thm}

\begin{rem}{\rm It is worth noticing that there is no hope for replacing the $\limsup$ in the last assertion of Theorem~\ref{thm-asympt-BS-dual-upper}~(\ref{assertion-minoration-dual}) by $\liminf$. Indeed, take $\psi \in H(\overline{\D})$ and set $\varphi(z)=\psi(z^d)$, for some integer $d\geq 2$. Then, for every $n\in \N$, $\widehat{\varphi^n}(k)=0$ whenever $d$ does not divide $k$, hence
	\[	
	\liminf_{k}(k+1)^{-p\alpha_p(\varphi)}\sum_{n\geq 0} |\widehat{\varphi^n}(k)|^p =0.
	\]
	In \cite{CharpentierEspoullierZarouf2026-1}, the authors establish a sharp version of Assertion~\ref{assertion-minoration-dual} for $\varphi$ a finite Blaschke product and $p=2$, with a self-contained proof, more elementary than the one displayed below.
	}
\end{rem}

Before proceeding to the proof of Theorem~\ref{thm-asympt-BS-dual-upper}, let us derive Theorem~\ref{thm-main-1-to-prove} from it.

\begin{proof}[Proof of Theorem~\ref{thm-main-1-to-prove}]If $p\geq 2$, one has
	\[
	\alpha_p(\varphi)\leq 0 \quad \text{and} \quad \sup_{k\geq 0}\sum _{n\geq 0} |\widehat{\varphi^n}(k)|^{p}<\infty,
	\]
	for every $\varphi\in H(\overline{\D})$. The first inequality is trivial by the definition of $\alpha_p(\varphi)$, while the second one follows from Proposition~\ref{prop-key-boundedness-carac}~(2) with $w=1$ and the fact that $C_\varphi:\,\ell_a^2\to \ell_a^2$ is bounded for every symbol $\varphi$.
	
The case $p\in [1,2]$ is a direct consequence of Assertion~\ref{eq-asympt-BS-sup-majoration} and the last assertion of Theorem~\ref{thm-asympt-BS-dual-upper}~(\ref{assertion-minoration-dual}).
\end{proof}

The remainder of the section is devoted to proving Theorem~\ref{thm-asympt-BS-dual-upper}. The proof of Assertion~\ref{eq-asympt-BS-sup-majoration} is somewhat lengthy and technical, and draws on several arguments developed by Blyudze and Shimorin in \cite{BlyudzeShimorin1993}. We shall use the following version of van der Corput's lemma.

\begin{lem}[Chapter VIII, Proposition~2 and the subsequent corollary in \cite{SteinBook1993}]\label{Stein-vanderCorput}
	Let $r\geq 1$, let $F\in C^r([a,b])$ be real-valued, and suppose that
	\[
	|F^{(r)}(t)|\geq \lambda>0,\quad a\leq t\leq b.
	\]
	If $r=1$, assume in addition that $F'$ is monotone on $[a,b]$.	Then, for every function $G\in C^1([a,b])$,
	\[
	\left| \int_a^b G(t)e^{iF(t)}\,dt \right| \leq C_r\lambda^{-1/r}
	\left(\sup_{t\in[a,b]}|G(t)| + \int_a^b |G'(t)|\,dt	\right),
	\]
	where $C_r>0$ depends only on $r$.
\end{lem}

Let us turn to the proof of Theorem~\ref{thm-asympt-BS-dual-upper}.

\begin{proof}[Proof of Theorem~\ref{thm-asympt-BS-dual-upper}]Assertion~\ref{assertion-minoration-dual} proceeds easily from Theorem~\ref{thm-asympt-BS}~(\ref{assert-BS-minoration}), so we shall focus on it first. Fix $p\in[1,\infty)$. By Theorem~\ref{thm-asympt-BS}~(\ref{assert-BS-minoration}) and Lemma~\ref{lem-elementary}, there exist $0<J<1<L<\infty$, $s\in(0,1)$ and $N_0\in \N$ such that, for every $n\geq N_0$,
\[
\sum_{Jn\leq k\le Ln}|\widehat{\varphi^n}(k)|^{p} + s^n \gtrsim \sum_{k\geq 0}|\widehat{\varphi^n}(k)|^{p}\gtrsim (n+1)^{p\alpha_{p}(\varphi)}.
\]
Upon enlarging $N_0$ if necessary, we get
\[
\sum_{Jn\leq k\le Ln}|\widehat{\varphi^n}(k)|^{p} \gtrsim (n+1)^{p\alpha_{p}(\varphi)},\quad n\geq N_0.
\]
Therefore, for every $N\geq N_0$,
\[
N^{1+p\alpha_{p}(\varphi)}\lesssim \sum_{n=N}^{2N} (n+1)^{p\alpha_{p}(\varphi)} \lesssim \sum_{JN\leq k\le 2LN}\sum_{n\geq 0} |\widehat{\varphi^n}(k)|^{p},
\]
which is Assertion~\ref{assertion-minoration-dual}. The last assertion of Theorem~\ref{thm-asympt-BS-dual-upper}~(\ref{assertion-minoration-dual}) is easily deduced by contradiction and the fact that
\[
\frac1N\sum_{JN\leq k\le 2LN}(k+1)^{p\alpha_p(\varphi)} \lesssim N^{p\alpha_p(\varphi)}.
\]

\medskip
	
Let us now deal with the proof of Assertion~\ref{eq-asympt-BS-sup-majoration}. Throughout the proof, the implicit constants in the notation $\lesssim$ and $\asymp$ are independent of $n$ and $k$. Let $p\in[1,2]$ be fixed. By Lemma~\ref{lem-elementary}, there exist $0<J<1<L<\infty$ and $s\in(0,1)$ such that
	\[
	\sum_{n\geq 0}|\widehat{\varphi^n}(k)|^p \lesssim \sum_{Jk \leq n\leq Lk}|\widehat{\varphi^n}(k)|^p +s^k,
	\]
	for every $k$ large enough. Thus, it is sufficient to prove that there exists $k_0\in \N$ such that, for every $k\geq k_0$,
	\begin{equation}\label{eq-estim-BS-dual-reduced}
		\sum_{Jk \leq n\leq Lk}\left|\widehat{\varphi^n}(k)\right|^p \lesssim	k^{p\alpha_p(\varphi)}.
	\end{equation}
	By Hölder's inequality, we get
	\[
	\sum_{Jk \leq n\leq Lk}|\widehat{\varphi^n}(k)|^p \lesssim k^{\frac{2-p}{2}} \left(\sum_{n\geq 0} |\widehat{\varphi^n}(k)|^2\right)^{\frac{p}{2}}.
	\]
	Then, denoting by $\left<\cdot,\cdot\right>$ the usual scalar product of $\ell_a^2$, the boundedness of $C_\varphi$ on $\ell_a^2$ yields
	\[
	\sum_{n\geq 0}|\widehat{\varphi^n}(k)|^2 = \sum_{n\geq 0}|\left<z^n,C_\varphi^*(z^k)\right>|^2=\Vert C_\varphi ^*(z^k)\Vert _{\ell_a^2} ^2 \leq \Vert C_\varphi\Vert ^2< \infty,
	\]
	hence
	\[
	\sum_{Jk \leq n\leq Lk}|\widehat{\varphi^n}(k)|^p \lesssim k^{\frac{2-p}{2}}.
	\]
	When $\varphi$ is a finite Blaschke product, the above estimate is precisely Assertion~\ref{eq-asympt-BS-sup-majoration}. When $\varphi$ is not a finite Blaschke product, the estimate Assertion~\ref{eq-asympt-BS-sup-majoration} is better than the above one and its proof requires much more work.
	
	\medskip
	
	In the remainder of the proof, we assume that $\varphi$ is not a finite Blaschke product. We shall first show that, in order to prove \eqref{eq-estim-BS-dual-reduced}, it is enough to replace $|\widehat{\varphi^n}(k)|$ by the contribution, in the integral defining the $k$th Fourier coefficient of $\varphi^n$, of the corresponding integral restricted to a small neighborhood of a given $\theta \in J_\varphi$. Recall that
	\[
	\widehat{\varphi^n}(k)=	\frac{1}{2\pi}
	\int_{-\pi}^{\pi}
	\varphi(e^{\mathrm{i}t})^n e^{-\mathrm{i}kt}\,dt.
	\]
	Since $\varphi$ is not a finite
	Blaschke product, Lemma~1.5 of \cite{BlyudzeShimorin1993} shows that
	$J_\varphi$ is finite. Write $J_\varphi=\{\theta_1,\ldots,\theta_N\}$.	Choose pairwise disjoint neighborhoods $U_1,\ldots,U_N$ of these
	points and functions $\chi_j\in C^\infty(\mathbb T)$ such that $0\leq\chi_j\leq1$, $\operatorname{supp}\chi_j\subset U_j$, and $\chi_j=1$ on a neighborhood of $\theta_j$. Set
	\[
	\chi_0=1-\sum_{j=1}^N\chi_j.
	\]
	Then
	\[
	\widehat{\varphi^n}(k)
	=
	\sum_{j=1}^N I_{n,k}^{(j)}+E_{n,k},
	\]
	where
	\[
	I_{n,k}^{(j)}=\frac{1}{2\pi}\int_{-\pi}^{\pi}	\chi_j(t)\varphi(e^{it})^n	e^{-ikt}\,dt
	\]
	and
	\[
	E_{n,k}	=\frac{1}{2\pi}	\int_{-\pi}^{\pi}	\chi_0(t)\varphi(e^{it})^n	e^{-ikt}\,dt.
	\]
	By definition of $J_\varphi$ and $\chi_0$, we have $|\varphi(e^{it})|<1$ for every $t\in\operatorname{supp}\chi_0$. The compactness of $\operatorname{supp}\chi_0$ and the continuity of
	$\varphi$ imply that
	\[
	s_0
	:=
	\max_{t\in\operatorname{supp}\chi_0}
	|\varphi(e^{it})|
	<1.
	\]
	Hence, for every $k$ and every $Jk\leq n\leq Lk$,
	\[
	\begin{aligned}
		|E_{n,k}|
		&\leq \frac{1}{2\pi}\int_{-\pi}^{\pi}|\chi_0(t)|\,|\varphi(e^{it})|^n\,dt \\
		&\leq s_0^n\\
		& \leq s^k,
	\end{aligned}
	\]
	where $s=s_0^J\in(0,1)$. Summing over $Jk\leq n\leq Lk$, we deduce
	\begin{align*}
		\sum_{Jk\leq n\leq Lk} 	|\widehat{\varphi^n}(k)|^p & \leq (N+1)^{p-1} \sum_{Jk\leq n\leq Lk}\left(\sum_{j=1}^N |I_{n,k}^{(j)}|^p + |E_{n,k}|^p\right) \\
		& \lesssim \sum_{j=1}^N \sum_{Jk\leq n\leq Lk} |I_{n,k}^{(j)}|^p	+ ks^k.
	\end{align*}
	
	For $\theta\in J_\varphi$, set
	\[
	\beta_p(\theta)=\frac{(\sigma(\theta)-r(\theta))(2-p)-2(p-1)}{2\sigma(\theta)}.
	\]
	By the previous estimate, in order to establish \eqref{eq-estim-BS-dual-reduced}, it therefore suffices to prove
	\begin{equation}\label{eq-main-thm-reduced}
		\sum_{Jk\leq n\leq Lk}
		|I_{n,k}^{(j)}|^p \lesssim k^{\beta_p(\theta_j)},\quad j \in \{1,\ldots, N\}.
	\end{equation}
	
	From now on, fix $\theta \in J_\varphi$ and let us simply write $r=r(\theta)$, $\sigma=\sigma(\theta)$, $\beta_p =\beta_p(\theta)$, $\chi(t)=\chi_j(t+\theta)$ and
	\[
	I_{n,k}=I_{n,k}^{(j)}=\frac{1}{2\pi}\int_{-\pi}^{\pi}	\chi(t-\theta)\varphi(e^{it})^n	e^{-ikt}\,dt,
	\]
	where $j \in \{1,\ldots,N\}$ is such that $\theta=\theta_j$.
	
	For $|t| < \eta$, $\eta>0$ sufficiently small, write
	\begin{equation*}\label{expansion-phi-proof}
		g(t):=\log\varphi(e^{i(\theta+t)})= c_0 + i\tau t+c_rt^r+c_{r+1}t^{r+1}+O(t^{r+2}),
	\end{equation*}
	where $c_0 =\log\varphi(e^{i\theta})$ satisfies $\Re{c_0}=0$. As observed in the proof of Lemma~2.3 in \cite{BlyudzeShimorin1993}, we have $\tau >0$.
	
	Moreover, since $|\varphi(e^{i(\theta +t)})| <1$ for every $t\in [-\eta,\eta]\setminus \{0\}$ and since $\sigma$ is the first index at which the real part of the coefficients of the above Taylor expansion does not vanish, we may and shall assume that $\eta>0$ is small enough so that
	\begin{equation}\label{eq:local-decay}
		\Re (g(t)-g(0))\leq -c|t|^\sigma,
		\quad |t|\le\eta,
	\end{equation}
	for some $c>0$ depending only on $g$ and $\eta$.
	
	 Without loss of generality, we may and shall assume that $\operatorname{supp}(\chi)\subset [-\eta ,\eta]$. Thus,
	\begin{equation*}\label{eq-reductionJnk}
		I_{n,k}=\frac{e^{ng(0)-ik\theta}}{2\pi}J_{n,k},
	\end{equation*}
	where
	\[
	J_{n,k}=e^{-ng(0)}\int_{\R}\chi(t)e^{ng(t)-ikt}\,dt.
	\]
	As the core of the proof, we shall establish
	\begin{equation}\label{eq:endpoints}
		\sum_{Jk\le n\le Lk}|J_{n,k}|
		\lesssim k^{\frac{\sigma-r}{2\sigma}} \quad \text{and}\quad
		\sum_{Jk\le n\le Lk}|J_{n,k}|^2
		\lesssim k^{-\frac{1}{\sigma}}.
	\end{equation}
	
	\medskip
	
	\subsection*{The case $r=\sigma$.} Write
	\[
	h(t):=g(t) - g(0)-i\tau t = c_\sigma t^\sigma+O(t^{\sigma+1}),\quad t\in[-\eta,\eta].
	\]
    Upon decreasing $\eta$ if necessary, there exists $C>0$ such that
	\begin{equation}\label{eq:local-decay-plus}
		|h(t)|\leq C|t|^\sigma,\quad
		|h'(t)| \leq C|t|^{\sigma-1},\quad
		|h''(t)|\leq C|t|^{\sigma-2},\quad t\in [-\eta,\eta].
	\end{equation}
	
	Now, setting $m=k-\tau n$, we have
	\[
	J_{n,k} = \int _{\R}\chi(t)e^{nh(t)}e^{-imt}\,dt,
	\]
	and, setting $l=k^{-1/\sigma}$, after the change of variables $t=ls$, we obtain
	\begin{equation}\label{eq:I-rescaled}
		J_{n,k} = l\int_{\mathbb R}H_{n,k}(s)e^{-imls}\,ds,
	\end{equation}
	where $H_{n,k}(s) =\chi(ls)e^{nh(ls)}$.
	
	We claim that
	\begin{equation}\label{eq:H-bounds}
		\int_{\mathbb R}|H_{n,k}(s)|\,ds + \int_{\mathbb R}|H_{n,k}''(s)|\,ds \lesssim 1.
	\end{equation}
	
	Indeed, by \eqref{eq:local-decay} and the fact that $\operatorname{supp}(\chi) \subset [-\eta,\eta]$, we first have
	\[
	|H_{n,k}(s)|=|\chi(ls)|e^{n\Re h(ls)} \lesssim e^{-cnl^\sigma|s|^\sigma}.
	\]
	Since $Jk\leq n\leq Lk$, we have $nl^\sigma\asymp1$. Thus $|H_{n,k}(s)| \lesssim e^{-c_1|s|^\sigma}$. Consequently,
	\[
	\int_{\mathbb R}|H_{n,k}(s)|\,ds \lesssim \int_{\mathbb R}e^{-c_1|s|^\sigma}\,ds \lesssim 1.
	\]
	
	We now estimate $\int_{\mathbb R}|H_{n,k}''(s)|\,ds$. Differentiating once gives
	\[
	H_{n,k}'(s)=l\chi'(ls)e^{nh(ls)}+nl\,\chi(ls)h'(ls)e^{nh(ls)}.
	\]
	Differentiating once more, we obtain
	\begin{multline*}
		H_{n,k}''(s) = l^2\chi''(ls)e^{nh(ls)}	+2nl^2\chi'(ls)h'(ls)e^{nh(ls)}	\\
		+nl^2\chi(ls)h''(ls)e^{nh(ls)}+n^2l^2\chi(ls)\bigl(h'(ls)\bigr)^2e^{nh(ls)}.
	\end{multline*}
	We estimate these four terms separately. Using $nl^\sigma \asymp 1$ and
	\[
	|h''(ls)|\lesssim |ls|^{\sigma-2} \quad \text{and} \quad |h'(ls)|^2\lesssim |ls|^{2\sigma-2},
	\]
	we obtain
	\[
	nl^2|\chi(ls)|\,|h''(ls)|e^{n\Re h(ls)}\lesssim  |s|^{\sigma-2}e^{-c|s|^\sigma},
	\]
	and
	\[
	n^2l^2|\chi(ls)|\,|h'(ls)|^2e^{n\Re h(ls)} \lesssim	|s|^{2\sigma-2}e^{-c|s|^\sigma}.
	\]
	The right hand-sides of both estimates are integrable on $\mathbb R$ since $\sigma\geq2$, hence
	\[
	nl^2\int_{\mathbb R}|\chi(ls)|\,|h''(ls)|e^{n\Re h(ls)}\,ds \lesssim 1
	\]
	and
	\[
	n^2l^2	\int_{\mathbb R}|\chi(ls)|\,|h'(ls)|^2e^{n\Re h(ls)}\,ds\lesssim 1.
	\]
	
	It remains to consider the terms involving derivatives of $\chi$. Let us define $\delta>0$ such that $\chi=1$ on $[-\delta,\delta]$. Then, $\chi'$ and $\chi''$ are supported on $\{t\in \R:\, \delta\leq |t|\leq \eta\}$. Thus, whenever $\chi'(ls)$ or $\chi''(ls)$ is nonzero, we have $|ls|\geq\delta$, and therefore
	\[
	e^{n\Re h(ls)}\leq e^{-cn\delta^\sigma} \leq e^{-c_2k},
	\]
	for some $c_2>0$, since $n\asymp k$. It yields
	\[
	l^2\int_{\mathbb R}|\chi''(ls)|e^{n\Re h(ls)}\,ds\lesssim e^{-c_2k} \lesssim 1.
	\]
	Similarly, since $h'$ is bounded on the support of $\chi'$,
	\[
	nl^2\int_{\mathbb R}|\chi'(ls)|\,|h'(ls)|e^{n\Re h(ls)}\,ds \lesssim e^{-c_2k} \lesssim 1.
	\]
	The last bound in both inequalities follows because $n\asymp k$. Combining the four	estimates proves \eqref{eq:H-bounds}.
	
	We now return to \eqref{eq:I-rescaled}. If $l|m|\leq1$, then
	\eqref{eq:H-bounds} immediately gives
	\[
	|J_{n,k}|\leq l\int_{\mathbb R}|H_{n,k}(s)|\,ds	\lesssim l.
	\]
	Suppose now that $l|m|>1$. Since $H_{n,k}$ is compactly supported, two integrations by parts yield
	\[
	\int_{\mathbb R}H_{n,k}(s)e^{-imls}\,ds	= -\frac{1}{m^2l^2}	\int_{\mathbb R}H_{n,k}''(s)e^{-imls}\,ds,
	\]
	hence
	\[
	J_{n,k}=-\frac{1}{m^2l}\int_{\mathbb R}H_{n,k}''(s)e^{-imls}\,ds,
	\]
	and, by \eqref{eq:H-bounds}, $|J_{n,k}|	\lesssim  l(l|m|)^{-2}$. Observing that $1\lesssim (1+l|m|)^{-2}$ if  $l|m|\leq1$ and $(l|m|)^{-2} \lesssim (1+l|m|)^{-2}$ if $l|m|>1$, the previous discussion gives
	\begin{equation*}\label{eq:r=sigma-pointwise}
		|J_{n,k}|\lesssim l(1+l|m|)^{-2}.
	\end{equation*}
	Recalling $m=k-\tau n$, we get
	\[
	\sum_{Jk\leq n \leq Lk} l(1+l|m|)^{-2}\lesssim 1 \quad \text{and} \quad \sum_{Jk\leq n \leq Lk} l^2(1+l|m|)^{-4}\lesssim l.
	\]
	Since $l=k^{-1/\sigma}$ and $r=\sigma$, these are precisely the two estimates in \eqref{eq:endpoints}.
	
	\medskip
	
	\subsection*{The case $r<\sigma$.} Set $\rho=k^{-1/r}$, $g - g(0)=u+iv$ with $u$ and $v$ real valued, and
	\[
	G_n(t)=\chi(t)e^{nu(t)},\quad t\in [-\eta,\eta].
	\]
	By \eqref{eq:local-decay}, $|G_n(t)|\leq 1$ for every $|t|\leq \eta$.
	Moreover, since
	\[
	G_n'(t)	=\chi'(t)e^{nu(t)}+n\chi(t)u'(t)e^{nu(t)},
	\]
	Since $\sigma$ is the first non-zero order in the Taylor expansion of $u$, Taylor's formula gives $|u'(t)|\lesssim |t|^{\sigma-1}$ for $|t|\leq\eta$. Together with \eqref{eq:local-decay}, this yields
	\[
	\int_{\mathbb R}|G_n'(t)|\,dt \lesssim 1+n\int_{-\eta}^\eta |t|^{\sigma-1}e^{-cn|t|^\sigma}\,dt	\lesssim1.
	\]
	
	On the other hand, let 	\[
	F(t):=nv(t)-kt,\quad |t|\leq \eta.
	\]
	By the definition of $r\geq 2$, we have $F^{(r)}(t)=nv^{(r)}(t)$,	and we may decrease $\eta$ so that
	\[
	|v^{(r)}(t)|\geq c_0>0, \quad |t|\leq\eta.
	\]
	It follows $|F^{(r)}(t)|\gtrsim  k$, since $n\asymp k$.
	
	Therefore, we can apply the Van der Corput Lemma~\ref{Stein-vanderCorput} and get
	\begin{eqnarray}\label{eq:vdc-rho}
		\label{defi-Jnk}|J_{n,k}| &= & \left| \int _{-\eta}^{\eta} G_n(t)e^{iF(t)}\,dt \right| \\
		\notag & \lesssim & k^{-1/r} \bigl( \sup_{t\in [-\eta,\eta]}|G_n(t)|+ \int_{-\eta}^{\eta}  |G_n'(t)|dt	\bigr)\\
		\notag & \lesssim & \rho.
	\end{eqnarray}
	
	Fix a fixed constant $A$, to be chosen sufficiently large, independently of $n$ and $k$. Since $m=k-\tau n$, there are $O(\rho^{-1})$ integers $n$ for which $|m|\le A/\rho$. By \eqref{eq:vdc-rho}, we get
	\begin{equation}\label{eq-step1}
		\sum_{\substack{Jk\leq n\le Lk \\ |m|\leq A/\rho}}|J_{n,k}|
		\lesssim 1 \quad \text{and}\quad
		\sum_{\substack{Jk\leq n\le Lk \\ |m|\leq A/\rho}}|J_{n,k}|^2
		\lesssim \rho.
	\end{equation}
	
	To finish the proof of \eqref{eq:endpoints}, it remains to estimate the integral
	\[
	\int _{-\eta}^{\eta} G_n(t)e^{iF(t)}\,dt
	\]
	in the case $|m|>A/\rho$. We shall do it in several steps. Define
	\begin{equation}\label{eq:T-lambda}
		T=\left(\frac{|m|}{k}\right)^{1/(r-1)}.
	\end{equation}
	In particular, $kT^{r-1}=|m|$. The definition of $m$, the fact that $Jk\leq n \leq Lk$ and the assumption $|m|>A/\rho$ yield
	\begin{equation*}\label{eq:T-rho}
		\rho \lesssim T  \leq T_0:= (1+\tau L)^{1/(r-1)}
	\end{equation*}
	where, as throughout the proof, the implicit constant is independent
	of $n$ and $k$. 
	
	Choose $\psi\in C_c^\infty(\mathbb R)$ such that
	\[
	0\leq\psi\leq1,\quad
	\psi=1 \ \text{on }[-c,c]
	\quad \text{and}\quad
	\operatorname{supp}\psi\subset[-2c,2c],
	\]
	where $c>0$ will be chosen sufficiently small, independently of $n$ and $k$. Since $T\leq T_0$,
	we may first choose $c$ so small that $2cT_0<\eta$. It follows that
	\[
	\operatorname{supp}\psi(\,\cdot\,/T)
	\subset[-2cT,2cT]
	\subset[-\eta,\eta],
	\]
	so that all the local estimates obtained above may be used on $\operatorname{supp}\psi(\,\cdot\,/T)$.
	
	From the local expansion of $v$ we can write
	\begin{equation}\label{eq:P-def}
		v'(t)-\tau=t^{r-1}P(t),
	\end{equation}
	where $P$ is real-analytic in $[-\eta,\eta]$ and $P(0)\neq0$. If $t\in \operatorname{supp}\psi(\,\cdot\,/T)$, then
	\[
	\begin{aligned}
		|n(v'(t)-\tau)|
		&= n |t|^{r-1}|P(t)| \\
		& \lesssim k(cT)^{r-1}\\
		& \lesssim c^{r-1}|m|
	\end{aligned}
	\]
	We now fix $c>0$, independently of $n$ and $k$, so small that, for every $t\in \operatorname{supp}\psi(\,\cdot\,/T)$,
	\begin{equation}\label{eq:Fprime-lower-inner}
		\begin{aligned}
			|F'(t)|
			&= |n(v'(t)-\tau)-m| \\
			&\geq |m|-|n(v'(t)-\tau)| \\
			&\geq \frac12|m|.
		\end{aligned}
	\end{equation}
	
	We shall also need estimates for $F''$ and $F'''$. By \eqref{eq:P-def}, we have
	\[
	v''(t)=(r-1)t^{r-2}P(t)+t^{r-1}P'(t),
	\]
	hence, for every $t\in \operatorname{supp}\psi(\,\cdot\,/T)$,
	\begin{equation}\label{eq:Fthird-inner-0}
		|F''(t)|=	n|v''(t)| \leq n|t|^{r-2}\lesssim kT^{r-2}=|m|T^{-1}.
	\end{equation}
	Similarly,
	\begin{equation}\label{eq:Fthird-inner}
		|F'''(t)|\lesssim |m|T^{-2}.
	\end{equation}
	Above and below, we recall that all the implicit constants may depend on $r$, $J$, $L$, $\varphi$, $\chi$, and $\psi$, but not on $n$ or $k$.
	
	For $t\in[-\eta,\eta]$, set
	\[
	b_T(t)=\psi(t/T)\chi(t) \quad \text{and} \quad a(t)=b_T(t)e^{nu(t)}.
	\]
	We shall next establish
	\begin{equation}\label{eq:a-derivatives}
		\int_{-\eta}^{\eta} |a(t)|\,dt \lesssim T, \quad \int_{-\eta}^{\eta} |a'(t)|\,dt\lesssim1 \quad \text{and} \quad \int_{-\eta}^{\eta} |a''(t)|\,dt\lesssim\rho^{-1}.
	\end{equation}
	
	By \eqref{eq:local-decay}, we first have $|a(t)|\leq|b_T(t)|$, $|t|\leq \eta$.	So the first estimate in \eqref{eq:a-derivatives} follows from the fact that $b_T$ is supported in $[-2cT,2cT]$. 
	
	For the two other estimates, let us differentiate $b_T$:
	\[
	b_T'(t)=T^{-1}\psi'(t/T)\chi(t)+\psi(t/T)\chi'(t).
	\]
	Then we have $\int_{-\eta}^{\eta}|b_T'|\lesssim 1$ and, since
	\[
	a'(t)=b_T'(t)e^{nu(t)}+nb_T(t)u'(t)e^{nu(t)},
	\]
	we deduce that
	\[
	\int_{-\eta}^{\eta}|a'(t)|\lesssim 1.
	\]
	For the third estimate in \eqref{eq:a-derivatives}, Taylor's formula also gives $|u''(t)|\lesssim |t|^{\sigma-2}$ for $|t|\leq\eta$, and we have
	\[
	a''(t)=b_T''(t)e^{nu(t)}+2n b_T'(t)u'(t)e^{nu(t)} +n b_T(t)u''(t)e^{nu(t)} + n^2b_T(t)(u'(t))^2e^{nu(t)}.
	\]
	First,
	\[b_T''(t)=T^{-2}\psi''(t/T)\chi(t)+2T^{-1}\psi'(t/T)\chi'(t) +\psi(t/T)\chi''(t).
	\]
	Hence, by changes of variables in the terms involving derivatives
	of $\psi$,
	\[
	\int_{-\eta}^{\eta}|b_T''(t)|\,dt \lesssim	T^{-1}+1.
	\]
	Since $T\gtrsim\rho$ and $\rho\leq1$, we obtain
	\begin{equation*}\label{eq:bsecond}
		\int_{-\eta}^{\eta}|b_T''(t)|e^{nu(t)}\,dt\lesssim\rho^{-1}.
	\end{equation*}
	
	Similarly, we check that
	\[
	\int_{-\eta}^{\eta}|n b_T(t)u''(t)e^{nu(t)}|\,dt + \int_{-\eta}^{\eta}|n^2b_T(t)(u'(t))^2e^{nu(t)}|\,dt \lesssim k^{1/\sigma} \leq k^{1/r} = \rho^{-1}.
	\]
	For the second term in $a''$, one may use
	\[
	b_T'(t)
	=
	T^{-1}\psi'(t/T)\chi(t)
	+
	\psi(t/T)\chi'(t)
	\]
	to also verify that
	\[
	\int_{-\eta}^{\eta}|2n b_T'(t)u'(t)e^{nu(t)}|\,dt \lesssim \rho^{-1}.
	\]
	Thus \eqref{eq:a-derivatives} is proved.
	
	Since $e^{iF}=(iF')^{-1}(e^{iF})'$, integration by parts twice gives
	\begin{align*}
		\int_{-\eta}^{\eta} a e^{iF}
		&= -\int_{-\eta}^{\eta}\left(\frac{a}{iF'}\right)'e^{iF}\\
		&= \int_{-\eta}^{\eta} \left[\frac{1}{iF'}
		\left(\frac{a}{iF'}\right)'\right]'e^{iF},
	\end{align*}
	with no boundary terms because $a$ is compactly supported. Expanding the derivatives and using \eqref{eq:Fprime-lower-inner}--\eqref{eq:Fthird-inner-0}--\eqref{eq:Fthird-inner}, we obtain
	\begin{equation*}\label{eq:ibp-expanded-bound}
		\left|\int_{-\eta}^{\eta} a(t)e^{iF(t)}\,dt\right|
		\lesssim
		\frac{\int_{-\eta}^{\eta} |a(t)|\,dt}{|m|^2T^2}
		+\frac{\int_{-\eta}^{\eta} |a'(t)|\,dt}{|m|^2T}
		+\frac{\int_{-\eta}^{\eta} |a''(t)|\,dt}{|m|^2}.
	\end{equation*}
	Using \eqref{eq:a-derivatives} and $T\gtrsim\rho$ yields
	\[
	\left|\int_{-\eta}^{\eta} a(t)e^{iF(t)}\,dt\right|	\lesssim \frac{\rho^{-1}}{|m|^2},
	\]
	hence, since $\rho|m|>A>1$,
	\begin{equation}\label{eq:inner}
		\left|\int_{-\eta}^{\eta}\psi(t/T)G_n(t)e^{iF(t)}\,dt\right|\lesssim \rho(1+\rho|m|)^{-2}.
	\end{equation}
	
	\medskip
	
	It remains to estimate
	\[
	\left|\int_{-\eta}^{\eta}(1-\psi(t/T))G_n(t)e^{iF(t)}\,dt\right|.
	\]
	To do so, set
	\[
	K(t)=(1-\psi(t/T))G_n(t)=(1-\psi(t/T))\chi(t)e^{nu(t)}.
	\]
	Note that
	\begin{equation}\label{supportG}
		\operatorname{supp}(K) \subset [-\eta,-cT]\cup[cT,\eta]\quad\text{and} \quad \operatorname{supp}(\psi'(\cdot/T)) \subset \{t\in \R:\, cT\leq |t|\leq 2cT\}.
	\end{equation}
	Therefore, we obtain, for some constants $c_2,c_3>0$ independent of $n$ and $k$,
	\begin{align}
		\sup_{t\in[-\eta,\eta]}|K(t)|+ \int_{-\eta}^{\eta}|K'(t)|dt
		&\lesssim e^{-c_2kT^\sigma}+n\int_{cT\le|t|\le\eta}|t|^{\sigma-1}e^{-c_2n|t|^\sigma}\,dt \notag\\
		&\lesssim e^{-c_3kT^\sigma}.
		\label{eq:G-var}
	\end{align}
	
	Now, set $\lambda=(kT^{r-2})^{-1/2}$. Observe that, by $|m|>Ak^{1/r}\geq c_0>0$,
	\[
	\bigl(n(1+|m|)^{r-2}\bigr)^{-1/(2(r-1))}\asymp k^{-1/(2(r-1))}|m|^{-(r-2)/(2(r-1))}=\lambda.
	\]
	Thus, since
	\[
	F'(t)=nt^{r-1}P(t)-m,
	\]
	we may apply Lemma~1.7 in \cite{BlyudzeShimorin1993} with $b=|m|$ (and with $P$ replaced by $-P$ if $m<0$), to get the existence of two intervals $E_1$ and $E_2$ such that
	\[
	E:=\{t\in [-\eta,\eta]:\, |F'(t)|\leq \lambda^{-1}\} \subset E_1 \cup E_2,
	\]
	and the lengths of both $E_1$ and $E_2$ are $\lesssim \lambda$.
	
	Then, by \eqref{eq:local-decay} and \eqref{eq:G-var},
	\begin{align}
		\left|\int_{E_1 \cup E_2} K(t)e^{iF(t)}\,dt\right|\notag
		&\leq \int_{E_1 \cup E_2} |K(t)|\,dt  \\\notag
		&\leq	|E_1 \cup E_2|\sup_{t\in{[-\eta,\eta]}}|K(t)|  \\
		&\lesssim \lambda e^{-c_3kT^\sigma}. \label{eq:exceptional}	
	\end{align}

	Let us now estimate
	\[
	\int_W K(t)e^{iF(t)}\,dt,
	\]
	where $W:=\operatorname{supp}K\setminus (E_1 \cup E_2)$.	Since $E\subset E_1 \cup E_2$,
	\begin{equation}\label{eq:Fprime-nonexceptional}
		|F'(t)|\ge \lambda^{-1},\quad t\in W.
	\end{equation}
	We shall now decompose $W$ into a number of intervals, uniformly bounded with respect to $n$ and $k$, on each of which $F'$ is monotone. We have
	\[
	F''(t)=nt^{r-2}\bigl((r-1)P(t)+tP'(t)\bigr).
	\]
	Since $(r-1)P(0)\neq 0$, by continuity, $F'$ is monotone on $[-\eta,0]$ and $[0,\eta]$, up to decreasing $\eta$.
	
	Moreover, by \eqref{supportG}, we can cover $W$ by pairwise disjoint intervals $W_1,\ldots,W_N$, with $N\leq 6$, on each of which $F'$ is monotone and, by \eqref{eq:Fprime-nonexceptional},
	\[
	\inf_{t\in W_j}|F'(t)|\geq \lambda^{-1}.
	\]
	Therefore, we may apply Lemma~\ref{Stein-vanderCorput} to $F$ and $K$  on each interval $W_j$, to get
	\[
	\left|\int_{W_j}K(t)e^{iF(t)}\,dt\right|\lesssim \lambda \left(\sup_{t\in[-\eta,\eta]}|K(t)|+ \int_{-\eta}^{\eta}|K'(t)|dt	\right).
	\]
	By \eqref{eq:G-var}, we thus obtain
	\begin{equation}\label{eq:nonexceptional}
		\left|\int_W K(t)e^{iF(t)}\,dt	\right|	\lesssim\lambda e^{-c_3kT^\sigma}.
	\end{equation}
	
	Combining \eqref{eq:exceptional} and
	\eqref{eq:nonexceptional}, we conclude that
	\[
	\left|	\int_{-\eta}^{\eta}K(t)e^{iF(t)}\,dt\right|	\lesssim \lambda e^{-c_3kT^\sigma}.
	\]
	
	Combining this with \eqref{eq:inner} and \eqref{defi-Jnk}, we arrive at the estimate
	\begin{equation*}\label{eq:outer-pointwise}
		|J_{n,k}|\lesssim \rho(1+\rho|m|)^{-2}+\lambda e^{-c_3kT^\sigma},\quad |m|>A/\rho.
	\end{equation*}
	
	First, since $m=k-\tau n$, we have
	\[
	\sum_{Jk\leq n\leq Lk}\rho(1+\rho|m|)^{-2}\lesssim1\quad\text{and}\quad	\sum_{Jk\leq n\leq Lk}\rho^2(1+\rho|m|)^{-4}\lesssim\rho.
	\]
	It remains to sum the second term. Put $M=k^{1-\frac{r-1}{\sigma}}$. A direct calculation from \eqref{eq:T-lambda} gives
	\[
		\lambda=k^{-\frac{1}{2(r-1)}}|m|^{-\frac{r-2}{2(r-1)}}\quad \text{and} \quad	kT^\sigma=\left(\frac{|m|}{M}\right)^{\frac{\sigma}{r-1}}.
	\]
	Since $m=k-\tau n$, comparison with an integral yields
	\begin{align*}
		\sum_{\substack{Jk\leq n\le Lk \\ |m|> A/\rho}}\lambda e^{-c_3kT^\sigma}
		&\lesssim k^{-\frac{1}{2(r-1)}}\int_0^\infty x^{-\frac{r-2}{2(r-1)}}e^{-c_3(x/M)^{\frac{\sigma}{r-1}}}\,dx\\
		&\lesssim k^{-\frac{1}{2(r-1)}}M^{\frac{r}{2(r-1)}}=k^{\frac{\sigma-r}{2\sigma}},
	\end{align*}
	and similarly
	\[
	\sum_{\substack{Jk\leq n\le Lk \\ |m|> A/\rho}}\lambda^2e^{-2c_3kT^\sigma}\lesssim k^{-\frac{1}{r-1}}M^{\frac{1}{r-1}}=k^{-\frac{1}{\sigma}}.
	\]
	Since $r<\sigma$, we have $\rho=k^{-1/r}\le k^{-1/\sigma}$. This yields
	\[
		\sum_{\substack{Jk\leq n\le Lk \\ |m|> A/\rho}}|J_{n,k}|\lesssim k^{\frac{\sigma-r}{2\sigma}} \quad \text{and}
		\quad\sum_{\substack{Jk\leq n\le Lk \\ |m|> A/\rho}}|J_{n,k}|^2\lesssim k^{-\frac{1}{\sigma}}.
	\]
	Together with \eqref{eq-step1}, this proves \eqref{eq:endpoints}.
	
	Finally, by Hölder's inequality, we get
	\begin{align*}
		\sum_{Jk\le n\leq Lk}|J_{n,k}|^p
		&=\sum_{Jk\leq n\leq Lk}|J_{n,k}|^{2-p}|J_{n,k}|^{2p-2}\\
		&\leq\left(\sum_{Jk\leq n\leq Lk}|J_{n,k}|\right)^{2-p}\left(\sum_{Jk\leq n\leq Lk}|J_{n,k}|^2\right)^{p-1}\\
		&\lesssim	k^{\frac{(\sigma-r)(2-p)-2(p-1)}{2\sigma}}\\
		& = k^{\beta_p}.
	\end{align*}
	This proves \eqref{eq-main-thm-reduced} and concludes the proof of the theorem.
\end{proof}

\section{Unboundedness of the matrix $A_\varphi$ on $\ell^p$ when $r(\theta)<\sigma(\theta)$ for some $\theta \in J_\varphi$}\label{sec-thmfull}

This section deals with the proof of Theorem~\ref{main-thm-rlesssigma}. In the sequel, for convenience, we shall consider the action on $\ell^p$ of the transpose of $A_\varphi$,
\[
\AA_\varphi:=A_\varphi ^t.
\]
Clearly, the boundedness of $\AA_\varphi$ on $\ell^p$ is equivalent to that of $C_\varphi$ on $\ell_a^p$. Recall that $C_\varphi$ is bounded on $\ell_a^2$. Therefore, Theorem~\ref{main-thm-rlesssigma} will follow once we have proven the following more precise statement. We shall call a finite rank coordinate projection of $\ell^p$ any map $P :\, \ell^p \to \ell^p$ of the form
\[
P((a_n)_n) = (b_n)_n \text{ with } \left\{
\begin{array}{lll}
	b_n = a_n & \text{if } n\in I \\
	b_n = 0 & \text{if } n\notin I, 
\end{array}
\right.
\]
where $I = [c,C]\cap \N$ for some $c,C>0$. We shall also write $\ell^p(\Z)$ for the space of sequences $u=(u_n)_{n\in \Z}$ such that
\[
\Vert u \Vert^p_{\ell^p(\Z)}= \sum_{n\in \Z} |u_n|^p <\infty.
\]
We keep on using the notation $\ell^p$ for $\ell^p(\N)$.

\begin{thm}\label{thm:local-matrix-obstruction}
	Let $\varphi\in H(\overline{\D})$ be a symbol such that $J_\varphi \neq \emptyset$. Assume that there exists $\theta_0\in J_\varphi$ such that
	\[
	r:=r(\theta_0)<\sigma:=\sigma(\theta_0)<\infty.
	\]
	Then, there exist sequences of finite rank coordinate projections $P_N$ of $\ell^p$ and of uniformly bounded convolution operators
	\[
	L_N:\ell^p\longrightarrow\ell^p(\Z),
	\]
	such that, for every $p\in[1,\infty]$ and every $N$ sufficiently large,
	\begin{equation*}\label{eq:quantitative-obstruction}
		\|L_N\mathcal A_\varphi P_N\|_{\ell^p\longrightarrow\ell^p(\Z)}\gtrsim N^{\frac{\sigma-r}{\sigma}|\frac{1}{p}-\frac{1}{2}|}.
	\end{equation*}
	In particular, $\mathcal A_\varphi$ is unbounded on $\ell^p$ for every $p\neq2$.
\end{thm}

The terminology ``convolution operators'' is made clear in the following elementary lemma. For $a<b$, we shall denote by $C^{\infty}_{c}((a,b))$ the space of $C^{\infty}$ functions with compact support in $(a,b)$. Given a sequence $u=(u_n)_{n\in\N}$, we shall also denote by $\tilde{u}$ the sequence indexed in $\Z$ and defined by
\[
\tilde{u}(n)=\left\{
\begin{array}{ll}
	u_n & \text{if }n \in \N\\
	0 & \text{if } n \in \Z\setminus \N
\end{array}.
\right.
\]

\begin{lem}\label{lem:uniform-localization} Let $a<b$, $\eta\in C^{\infty}_{c}((a,b))$, $\theta_{0}\in[0,2\pi)$. Let also $0<T\leq 1$ be such that $(b-a)T<2\pi$, and set
\[
\mathcal{U}=\{(\theta_{0}+t)\bmod2\pi:\ aT<t<bT\}\subset[0,2\pi).
\]
 Let $\eta_T$ be the $2\pi$-periodic function on $\R$ defined by 
 \begin{equation}\label{defi-etaN}
  \eta_T(\theta)=\left\{
  \begin{array}{ll}
\eta(t/T) & \text{if }\theta \in \mathcal{U}\\
0 & \text{if }\theta \in [0,2\pi) \setminus \mathcal{U}
 \end{array}
 \right.,
 \end{equation}
 where $t$ is the unique element of $(aT,bT)$ such that $\theta\equiv\theta_{0}+t\pmod{2\pi}$. Then, the convolution operator $L_T:\ell^{p}\longrightarrow\ell^{p}(\Z)$, defined by
\begin{equation}\label{defi-conv-op-L_T}
L_T(u)= \widehat{\eta_T} * \tilde{u},\quad \widehat{\eta_T} * \tilde{u}(k)=\sum_{m\in\Z}\widehat{\eta_T}(m)\tilde{u}(k-m),
\end{equation}
is bounded and satisfies
\[
\|L_T(u)\|_{\ell^{p}(\Z)}\leq C_{\eta}\|u\|_{\ell^{p}},
\]
where the constant $C_{\eta}$ depends only on $\eta$, and not on $T$ and $p\in[1,\infty]$.
 \end{lem}

\begin{proof}Since $\eta$ is compactly supported in $(a,b)$, the function $\eta_T$ is $C^{\infty}$ on $\R$. For $k\in\Z$, by definition and after the change of variable $t=Ts$, we have 
\[
\widehat{\eta_T}(k)=\frac{1}{2\pi}\int^{bT}_{aT}\eta(t/T)e^{-ik(\theta_{0}+t)}\,dt=\frac{Te^{-ik\theta_{0}}}{2\pi}\int^{b}_{a}\eta(s)e^{-ikTs}\,ds.
\]
It follows
\[
|\widehat{\eta_T}(k)|\le\frac{T}{2\pi}\|\eta\|_{L^{1}(a,b)}.
\]
Moreover, for $m\neq 0$, let $\omega=mT$. Using again that $\eta$ is compactly supported in $(a,b)$, two integrations by parts yield
\[
\int^{b}_{a}\eta(s)e^{-i\omega s}\,ds=\frac{1}{(i\omega)^{2}}\int^{b}_{a}\eta''(s)e^{-i\omega s}\,ds.
\]
 Thus, for every $m\neq 0$, 
\[
|\widehat{\eta_T}(k)|\le\frac{T}{2\pi}\frac{\|\eta''\|_{L^{1}(a,b)}}{(|k|T)^{2}}.
\]
 The two previous inequalities implies the existence of a constant $K$, depending only on $\eta$, such that
\[
|\widehat{\eta_T}(k)|\le KT(1+T|k|)^{-2},
\]
hence
\[
\Vert \widehat{\eta_T}\Vert_{\ell^{1}(\Z)}\leq KT\sum_{k\in\Z}(1+T|k|)^{-2}\le KT+2KT\int^{\infty}_{0}(1+Tx)^{-2}\,dx=KT+2K\leq C_\eta,
\]
where $C_\eta$ depends only on $\eta$. By Young's inequality, we get
\[
\|L_T(u)\|_{\ell^{p}(\Z)}=\|\widehat{\eta_T}*\tilde{u}\|_{\ell^{p}(\Z)}\le\|\widehat{\eta_T}\|_{\ell^{1}(\Z)}\|\tilde{u}\|_{\ell^{p}(\Z)}\le C_{\eta}\|u\|_{\ell^{p}(\N)},
\]
 as desired.

\end{proof}

The main ingredient of the proof of Theorem~\ref{thm:local-matrix-obstruction} is a uniform version of the stationary phase method, due to Fedoryuk. The next statement is a reformulation of \cite[Theorem~1.6, p.~107]{Fedoryuk1977}, that can also be found in \cite[Theorem~4]{AuxemeryBorichevZarouf2026}.

\begin{thm}[Theorem~1.6, p.~107 of \cite{Fedoryuk1977}]\label{thm-Feodryuk}For $d\geq 1$, let $\Omega \subset \R^d$ be an open set. Let also $a<b$ and let
	\[
	A:\,[a,b]\times \Omega \to \C\quad \text{and}\quad \Phi:\,[a,b]\times \Omega \to \R
	\]
be two functions satisfying the following hypotheses:
\begin{enumerate}[(i)]
	\item \label{item1Fedor} $A$ is $C^{\infty}$ on $[a,b]\times \Omega$ and, together with all its derivatives in the first variable, vanishes at $(a,\alpha)$ and $(b,\alpha)$ for every $\alpha \in \Omega$;
	\item \label{item2Fedor} $\Phi$ is $C^\infty$ on $[a,b]\times \Omega$ and, for every $\alpha \in \Omega$, the function $\Phi_{\alpha}(s):=\Phi(s,\alpha)$ is real-valued and possesses a unique non-degenerate stationary point $s_*(\alpha) \in(a,b)$ satisfying, for some $\delta_0>0$,
	\[
	|\Phi_{\alpha}''(s_*(\alpha))|\geq \delta_0,\quad \alpha \in \Omega.
	\]
	Then, as $n\to \infty$,
	\begin{multline*}
		\int^{b}_{a}A(s,\alpha)e^{in\Phi_{\alpha}(s)}\,ds\\=\left(\frac{2\pi}{n|\Phi_\alpha''(s_*(\alpha))|}\right)^{1/2}A(s_{*}(\alpha),\alpha)e^{in\Phi_{\alpha}(s_*(\alpha))+i\frac{\pi}{4}\operatorname{sgn}\Phi_\alpha''(s_*(\alpha))}+O(n^{-3/2}),
	\end{multline*}
	where the constant in the $O(n^{-3/2})$ is uniform over $\alpha$ in all compact subsets of $\Omega$.
	
\end{enumerate}

\end{thm}

We now turn to the proof of Theorem~\ref{thm:local-matrix-obstruction}.

\begin{proof}[Proof of Theorem~\ref{thm:local-matrix-obstruction}]
	
Let $0<a<1<b$ and let $\eta \in C_c^{\infty}((a,b))$ be equal to $1$ on an open neighborhood of $1$. For a sequence $(T_N)_{N\geq 0}$ of positive real numbers such that $T_N\to 0$ as $N\to \infty$, let us denote by $\eta_N:= \eta_{T_N}$ the functions defined by \eqref{defi-etaN}, and let $L_N:=L_{T_N}$ be the convolution operators defined as in \eqref{defi-conv-op-L_T}
	
\medskip

We will prove that there exist such a sequence $(T_N)_{N\geq 0}$ such that
\begin{equation}\label{todo}
\|L_N\mathcal A_\varphi\|_{\ell^p\to\ell^p(\Z)}\gtrsim N^{\frac{\sigma-r}{\sigma}|\frac{1}{p}-\frac{1}{2}|},
\end{equation}
for every $N$ sufficiently large. The strategy of the proof will provide a sequence of finite rank projections $P_N$, as in the statement. By Lemma~\ref{lem:uniform-localization}, the last assertion of the statement would immediately follow.

\medskip

Let $\theta_{0}\in J_\varphi$ be such that $r:=r(\theta_0)<\sigma=\sigma(\theta_0)$, as in the statement of the theorem. For reasons that will make the sequel technically slightly easier, set
\[
T_N = N^{-1/\sigma}.
\]
Fix $N_0\in \N$ such that $(b-a)T_{N_0} < 2\pi$. This properly defines the $C^{\infty}$ functions $\eta_N$, for every $N\geq N_0$, and the associated convolution operators $L_N$. We will prove that this choice works for our purpose.

\medskip

In order to define the sequence $(u^N)_{N\geq 0}$ suitably, we must first understand the asymptotic behaviour of the coefficients of the matrix $L_N\mathcal A_\varphi$. We shall write them as oscillatory integrals, which we estimate using Theorem ~\ref{thm-Feodryuk}. Now, by definition of $L_N$, it is easily checked that
\[
L_N\mathcal A_\varphi({k,n})=\widehat{\eta_{N}\varphi^{n}}(k), \quad k\in \Z,n \in \N,
\]
where $L_N\mathcal A_\varphi({k,n})$ stands for the coefficients at line $k$ and column $n$ of $L_N\AA_\varphi$.

\medskip

To get a relevant expression of $\widehat{\eta_{N}\varphi^{n}}(k)$, let us write, as in the proof of Theorem~\ref{thm-asympt-BS},
\begin{equation*}\label{expansion-phi-proofbis}
	\log\varphi(e^{i(\theta_0+t)})= c_0 + i\tau t+u(t)+iv(t),\quad |t| < \rho,
\end{equation*}
where $\rho>0$ is sufficiently small so that $u,v$ are real analytic in a neighborhood of $[-\rho,\rho]$, and where $c_0 =\log\varphi(e^{i\theta_0})$. We recall that $\tau >0$. By definition of $r<\sigma$, we have
\[
u(t)=t^{\sigma}\widetilde{u}(t)\quad \text{and} \quad v(t)=t^{r}\widetilde{v}(t), \quad t\in[-\rho,\rho]
\]
with $\widetilde{u}(0)<0$ and $\sigma$ even (since $\varphi(\D)\subset \D$), and $\widetilde{v}(0)\neq 0$. Upon enlarging $N_0$ if necessary, one may assume that $[aT_N,bT_N]\subset [-\rho,\rho]$. Therefore,
\begin{equation}\label{eq:bloc-integrale-fr}
	\begin{aligned}\widehat{\eta_{N}\varphi^{n}}(k) & =\frac{\varphi^n(e^{i\theta_0})e^{-ik\theta_{0}}}{2\pi}\int^{bT_{N}}_{aT_{N}}\eta(t/T_{N})e^{nu(t)}e^{i[nv(t)-(k-\tau n)t]}\,dt\\
		& = \frac{T_{N}\varphi^n(e^{i\theta_0})e^{-ik\theta_{0}}}{2\pi}\int^{b}_{a}\eta(s)e^{\frac{n}{N}s^{\sigma}\widetilde{u}(T_N s)}e^{i\Lambda_N\Phi_{\alpha_{N,n,k}}(s)}\,ds,
	\end{aligned}
\end{equation}
where
\begin{equation}\label{eq:alpha-fr}
\Lambda_N=NT_N^r=N^{1-r/\sigma}\quad \text{and} \quad\alpha_{N,n,k}=\left(\frac{n}{N},\frac{k-\tau n}{NT_N^{r-1}},T_{N}\right),
\end{equation}
and
\[
\Phi_{\alpha}(s)=\mu s^{r}\widetilde{v}(\varepsilon s)-\nu s,\quad \alpha=(\mu,\nu,\varepsilon)\in\R^2\times (-1,1).
\]
Now, let $M_N$, $N\geq N_0$, be the matrix with coefficients given by
\begin{equation}\label{defi-MN}
M_N(k,n)=\frac{T_N}{2\pi}\int^{b}_{a}\eta(s)e^{\frac{n}{N}s^{\sigma}\widetilde{u}(T_N s)}e^{i\Lambda_N\Phi_{\alpha_{N,n,k}}(s)}\,ds.
\end{equation}
Since the matrices $U=\operatorname{diag}(\varphi^n(e^{i\theta_0}))_{n\geq 0}$ and $V=\operatorname{diag}(e^{-ik\theta_{0}})_{k\in \Z}$ are isometries of $\ell^p$ and $\ell^p(\Z)$ respectively, and since $L_N\mathcal{A}_\varphi=VM_{N}U$, it is enough to prove \eqref{todo} with $L_N\mathcal{A}_\varphi$ replaced by $M_N$. From now on, we will thus focus only on $M_N$.

\medskip

Observe that $\Lambda_N \to +\infty$ as $N\to \infty$. In order to estimate $M_N(k,n)$, we shall first seek for an open set $\Omega \subset \R^2 \times(-1,1)$ such that we may apply Fedoryuk's result to the oscillatory integral
\[
\int^{b}_{a}A(s,\alpha)e^{in\Phi_{\alpha}(s)}\,ds,
\]
where $A : [a,b]\times \Omega \to \C$ is the amplitude function defined by
\[
A(s,\alpha) = \eta(s)e^{\mu s^\sigma\tilde{u}(\varepsilon s)},\quad \alpha =(\mu,\nu,\varepsilon)\in \Omega,
\]
and $\Phi : [a,b]\times \Omega \to \R$ the phase function defined by
\[
\Phi(s,\alpha)=\Phi_{\alpha}(s),\quad \alpha =(\mu,\nu,\varepsilon)\in \Omega.
\]
Note that the hypothesis \eqref{item1Fedor} of Theorem~\ref{thm-Feodryuk} is satisfied since $\eta$ is $C^\infty$ on $\R$ and compactly supported in $(a,b)$. To check the second hypothesis, we apply the implicit function theorem to $\partial_s\Phi = \Phi'_\alpha$. Set $\alpha_0=(1,r\widetilde{v}(0),0)$ and observe that
\begin{equation}\label{derivPhi}
\partial _s\Phi(s,\alpha_{0})=r\widetilde{v}(0)(s^{r-1}-1) \quad \text{and} \quad \partial^2_{s}\Phi(s,\alpha_{0})=r(r-1)\widetilde{v}(0)s^{r-2}.
\end{equation}
Then $\partial _s\Phi(s,\alpha_{0})$ vanishes at $s_*(\alpha_0)=1$. Moreover, there exists $c_0>0$ such that $|\partial^2_{s}\Phi(s,\alpha_{0})| \geq 2c_0$ for every $s\in[a,b]$ hence, by continuity of $(s,\alpha)\mapsto\partial^2_{s}\Phi(s,\alpha)$ and compactness of $[a,b]\times\{\alpha_{0}\}$, one has
\begin{equation}\label{derivsec}
|\partial^2_{s}\Phi(s,\alpha)| \geq c_0
\end{equation}
for every $(s,\alpha)\in [a,b]\times \Omega_0$ and some open neighborhood $\Omega_0$ of $\alpha_0$. Now \eqref{derivsec} and the implicit function theorem ensure that there exists an open set $\Omega \subset \Omega_0$ containing $\alpha_0$ such that, for every $\alpha \in \Omega$, the function $\partial _s\Phi(s,\alpha)$ vanishes at a unique point $s_*(\alpha) \in [a,b]$, and the function $\alpha \mapsto s_*(\alpha)$ is $C^\infty$ on $\Omega$. This shows that the hypothesis \eqref{item2Fedor} of Theorem~\ref{thm-Feodryuk} holds, and thus that, as $n\to +\infty$, one has
\begin{equation}\label{estim-Fedor}
	\int^{b}_{a}A(s,\alpha)e^{in\Phi_{\alpha}(s)}\,ds=\frac{\beta(\alpha)}{\sqrt{n}}e^{inS(\alpha)}+O(n^{-3/2}),
\end{equation}
where
\[
\beta(\alpha)=\left(\frac{2\pi}{|\Phi_\alpha''(s_*(\alpha))|}\right)^{1/2}A(s_{*}(\alpha),\alpha)e^{\frac{i\pi}{4}\operatorname{sgn}\Phi_\alpha''(s_*(\alpha))}\quad\text{and}\quad S(\alpha)=\Phi_{\alpha}(s_*(\alpha)).
\]
\medskip

The next step is to replace $\alpha$ in the latter expansion by
\[
\alpha_{N,n,k}=\left(\frac{n}{N},\frac{k-\tau n}{NT_N^{r-1}},T_{N}\right),
\]
in order to obtain an asymptotic expansion of \eqref{defi-MN}. We recall that $\Lambda_N\to +\infty$ as $N\to \infty$, thus we only need to determine those $N$, $k$ and $n$ for which $\alpha_{N,n,k}\in \Omega$. Let us choose $d>0$ small enough so that the cube
\[
\mathcal{P}=[1-d,1+d]\times[r\widetilde{v}(0)-d,r\widetilde{v}(0)+d]\times[-d,d]
\]
is contained in $\Omega$. Upon reducing $d$, we may assume that $s_*(\alpha)$ is contained in $\{s\in [a,b]:\,\eta(s)=1\}$ for every $\alpha \in \PP$ (recall that $s_*(\alpha_0)=1$ and $\eta =1$ in an open neighborhood of $1$). Now, set
\[
\delta=\frac{d}{2(1+\tau)}\quad \text{and}\quad K_{N}=\tau N+r\widetilde{v}(0)NT_N^{r-1},
\]
and let
\[
I_{N}=\{n\in\N:|n-N|\le\delta NT_N^{r-1}\}\quad \text{and} \quad J_{N}=\{k\in\Z:|k-K_{N}|\leq\delta NT_N^{r-1}\}.
\]
Note that, for every $(n,k)\in I_{N}\times J_{N}$,
\[
\left|\frac{n}{N}-1\right| \leq T_N^{r-1}
\]
and
\[
\left|\frac{k-\tau n}{NT_N^{r-1}}-r\widetilde{v}(0)\right|\leq \frac{|k-K_{N}|+\tau|n-N|}{NT_N^{r-1}}\le(1+\tau)\delta\leq d/2.
\]
One may thus enlarge $N_0$ if necessary to get $\alpha_{N,n,k}\in \PP$ for every $N\geq N_0$ and every $(n,k)\in I_N\times J_N$. Therefore, the expansion \eqref{estim-Fedor} yields
\begin{eqnarray*}
M_N(k,n) & = & \frac{T_N}{2\pi}\left(\frac{\beta(\alpha_{N,n,k})}{\sqrt{\Lambda_N}}e^{i\Lambda_NS(\alpha_{N,n,k})}+O(\Lambda_N^{-3/2})\right)\\
& = & \frac{T_N}{2\pi\sqrt{\Lambda_N}}\left(\beta(\alpha_{N,n,k})e^{i\Lambda_NS(\alpha_{N,n,k})} +R_N(k,n)\right),
\end{eqnarray*}
where
\[
|R_N(k,n)|\leq \frac{C}{\Lambda_N},
\]
for some constant $C>0$ independent of $(n,k)\in I_N\times J_N$.

\medskip

The last step consists in choosing suitably the finitely supported sequence $u^N=u^N(n)$ so that
\begin{equation}\label{todo2}
	\frac{\|M_N (u^N)\|_{\ell^p(\Z)}}{\|u^N\|_{\ell^p}}\gtrsim N^{\frac{\sigma-r}{\sigma}|\frac{1}{p}-\frac{1}{2}|},
\end{equation}
for every $N$ sufficiently large. Observe first that every $N\geq N_0$ and every $(n,k)\in I_N\times J_N$, one has $|\beta(\alpha_{N,n,k})|\geq c_\beta$ for some constant $c_\beta>0$, and
\[
\beta(\alpha_{N,n,k})=e^{i\gamma}|\beta(\alpha_{N,n,k})|, \quad \gamma=\frac{\pi}{4}\rm{sgn}(\widetilde{v}(0)).
\]
The first part of this assertion follows from the compactness of $\PP$, the continuity of $\beta$ on $\Omega$, and the fact that, for every $\alpha \in \PP$,
\[
A(s_*(\alpha),\alpha) = e^{\mu s^\sigma\tilde{u}(\varepsilon s)} >0.
\]
The second part follows from \eqref{derivPhi}, \eqref{derivsec} and the continuity of $\partial ^2_s \Phi(s,\alpha)$ on $[a,b]\times \Omega$.

We distinguish the cases $p\in [2,\infty]$ and $p\in[1,2)$. Let us first consider $p\in [2,\infty]$. Let $k_0\in \Z$ be the closest integer to $K_N$. Set
\begin{equation*}
u^N(n)=
	\begin{cases}e^{-i\gamma - i\Lambda_NS(\alpha_{N,n,k_0})} & \text{if }n\in I_N\\
		0 & \text{if }n\notin I_N.
	\end{cases}
\end{equation*}
One has $\Vert u^N \Vert _{\ell^p} = |I_N|^{1/p}$ if $p<\infty$, where $|I_N|$ denotes the cardinal of $I_N$. If $p=\infty$, then $\Vert u^N\Vert_{\ell^{\infty}} =1$. We have
\begin{equation}\label{lasteqpgeq1}
 (M_N(u^N))_k  = \frac{T_N}{2\pi\sqrt{\Lambda_N}}\left(\sum_{n\in I_N}(|\beta(\alpha_{N,n,k})|e^{i\Lambda_N(S(\alpha_{N,n,k})-S(\alpha_{N,n,k_0}))} + R_N(k,n)u^N(n))\right),
\end{equation}
By definition, $\Lambda_NS(\alpha_{N,n,k})$ and $\Lambda_NS(\alpha_{N,n,k_0})$ are both the minimum or both the maximum on $[a,b]$ (depending only on the sign of $\widetilde{v}(0)$) of  the functions
\[
s\mapsto nv(T_Ns) -(k-\tau n)T_Ns) \quad \text{and} \quad  s\mapsto nv(T_Ns) -(k_0-\tau n)T_Ns),
\]
respectively. Since
\[
nv(T_Ns) -(k-\tau n)T_Ns - (nv(T_Ns) -(k_0-\tau n)T_Ns) = T_Ns(k_0-k),
\]
we thus get
\begin{equation*}
	 |\Lambda_N(S(\alpha_{N,n,k})-S(\alpha_{N,n,k_0}))| \leq bT_N|k-k_0|.
\end{equation*}
Let now
\[
J_N'=\{k\in \Z: |k-k_0|\leq (2bT_N)^{-1}\}.
\]
Note that $J_N'\subset J_N$ for $N$ large enough. Using $|e^{it} -1|\leq |t|$, we thus obtain, for $k\in J'_N$,
\begin{equation*}\label{lasteqpgeq2}
	|(M_N(u^N))_k| \geq  \frac{T_N|I_N|}{2\pi\sqrt{\Lambda_N}}\left(\frac{c_\beta}{2} - \frac{C}{\Lambda_N}\right).
\end{equation*}
Thus, since $|I_N|\asymp NT_N^{r-1}$, for $k\in J_N'$ and $N$ sufficiently large,
\[
|(M_N(u^N))_k|  \gtrsim \Lambda_N^{1/2}.
\]
Further, $|J_N'| \asymp T_N^{-1}$, hence, for $N$ sufficiently large,
\begin{equation*}
\Vert M_N(u^N)\Vert _{\ell^p(\Z)}  \gtrsim  |J_N'|^{1/p}\Lambda_N^{1/2} = T_N^{-1/p}\Lambda_N^{1/2}.
\end{equation*}
We thus obtain, for $N$ sufficiently large,
\[
\frac{\Vert M_N(u^N)\Vert _{\ell^p(\Z)}}{\Vert u^N\Vert_{\ell^p}} \gtrsim \Lambda_N^{1/2-1/p}.
\]

It remains to consider the case $p\in[1,2)$. We proceed by duality. Let
\[
F_N=(M_N(k,n))_{k\in J_N,n\in I_N}.
\]
and let $F_N^*$ be the adjoint of $F_N$. Denote $p'\in (2,\infty]$ the conjugate exponent of $p$. It is enough to prove that
\[
\Vert F_N^* \Vert _{\ell^{p'}(J_N)\to \ell^{p'}(I_N)} = \Vert F_N \Vert _{\ell^{p}(I_N)\to \ell^{p}(J_N)} \gtrsim \Lambda_N^{1/p-1/2}.
\]
Let
\[
v^N(k)=
	\begin{cases} e^{i\gamma + i\Lambda_NS(\alpha_{N,N,k})} & \text{if }k\in J_N\\
	0 & \text{if }k\notin J_N
	\end{cases}
\]
One has $\Vert v^N \Vert _{\ell^{p'}(\Z)} = |J_N|^{1/p'}$ if $p'$ is finite, and $\Vert v^N \Vert _{\ell^{\infty}(\Z)} = 1$. As previously, for $n\in I_N$, we obtain
\begin{equation}\label{lasteqpleq3}
|(F_N^*(v^N))_n| \geq \frac{T_N}{2\pi\sqrt{\Lambda_N}}\left(\sum_{k\in J_N}|\beta(\alpha_{N,n,k})|\left(1- |\Lambda_N(S(\alpha_{N,n,k})-S(\alpha_{N,N,k}))|\right) - \frac{C|J_N|}{\Lambda_N}\right).
\end{equation}
Observing that, for $s\in [a,b]$,
\[
nv(T_Ns) - (k-\tau n)T_Ns - (Nv(T_Ns) - (k-\tau N)T_Ns) = (n-N)(v(T_Ns) +\tau T_Ns),
\]
with $v(t)=O(t^r)$, using that $T_N\leq 1$, and proceeding as above, we get
\[
\Lambda_N|S(\alpha_{N,n,k})- S(\alpha_{N,N,k})|\leq CT_N|n-N|,
\]
for some constant $C>0$. Let
\[
I_N'=\{n\in \N:\, |n-N| \leq (2CT_N)^{-1}\}.
\]
Note that $I_N'\subset I_N$ for $N$ large enough, and $|I'_N|\asymp T_N^{-1}$.  Thus, by \eqref{lasteqpleq3}, for $n\in I'_N$,
\[
|(F_N^*(v^N))_n|  \gtrsim \frac{T_N|J_N|}{\sqrt{\Lambda_N}}\left(\frac{1}{2} -  \frac{C}{\Lambda_N}\right),
\]
and, for $N$ sufficiently large,
\[
|(F_N^*(v^N))_n|  \gtrsim \Lambda_N^{1/2}.
\]
Therefore,
\[
\Vert F_N^*(v^N) \Vert _{\ell^{p'}(I_N)}  \gtrsim |I'_N|^{1/p'}\Lambda_N^{1/2} \gtrsim T_N^{-1/p'}\Lambda_N^{1/2}.
\]
We deduce that, for $N$ large enough,
\[
\frac{\Vert F_N^*(v^N) \Vert _{\ell^{p'}(I_N)} }{\Vert v^N \Vert _{\ell^{p'}(J_N)}} \gtrsim \Lambda_N^{1/2-1/p'}=\Lambda_N^{1/p-1/2}.
\]

\medskip

To conclude the proof, it suffices to define $P_N$ as the coordinate projection onto $I_N$.
 \end{proof}
 
\medskip

We are now ready to prove Theorem~ \ref{Thm-Fulldescription}. The case of non-monomial finite Blaschke products has already been settled after Corollary~\ref{cor-BF-intro}, so we consider only the remaining symbols. For monomial symbols or those for which $J_\varphi = \emptyset$, the proof is almost trivial. For the proof of the remaining cases, we need essentially only to adapt that of Theorem~\ref{thm:local-matrix-obstruction}, highlighting how flexible it is. For convenience, let us restate Theorem~ \ref{Thm-Fulldescription}.

\begin{thm}\label{thm:mixed-exponents}Let $\varphi\in H(\overline{\D})$ be a symbol. Then the following assertions hold.
	\begin{enumerate}
		\item If $J_\varphi=\varnothing$, then $C_\varphi:\ell_a^p\to\ell_a^q$ is bounded for every $p,q\in[1,\infty]$.
		\item If $\varphi = \mu z^n$ is a monomial, for some $|\mu|\leq 1$ and some $n\geq 1$, then
		\begin{enumerate}
			\item If $|\mu|< 1$, then $C_\varphi:\ell_a^p\to\ell_a^q$ is bounded for every $p,q\in[1,\infty]$;
			\item If $|\mu|=1$, then, for every $p,q\in[1,\infty]$, $C_\varphi:\ell_a^p\to\ell_a^q$ is bounded if and only if $p\leq q$.
		\end{enumerate}
		\item If $\varphi$ is a non-monomial finite Blaschke product, then $C_\varphi:\ell_a^p\to\ell_a^q$ is bounded if and only if $1\leq p \leq 2 \leq q \leq \infty$.
		\item Assume that $J_\varphi\neq\varnothing$ and that $\varphi$ is neither a monomial nor a finite Blaschke product. Set
		\[
		D_\varphi:=\max_{\theta\in J_\varphi}	\left(\sigma(\theta)-r(\theta)\right).
		\]
		Then, for every $p,q\in[1,\infty]$, the following assertions are equivalent:
		\begin{enumerate}
			\item $C_\varphi:\ell_a^p\to\ell_a^q$ is bounded;
			\item Both inequalities
			\[
			(D_\varphi+1)\frac1p-\frac1q -\frac{D_\varphi}{2}\geq 0 \quad \text{and}\quad \frac1p-(D_\varphi+1)\frac1q+ \frac{D_\varphi}{2}\geq 0
			\]
			hold, with the convention $1/\infty=0$.
		\end{enumerate}
	\end{enumerate}
\end{thm}

\begin{rem}{\rm Observe that Assertion~(4) implies that if $C_\varphi:\ell_a^p\to \ell_a^q$ is bounded, then $p\leq q$. Moreover, if $D_\varphi=0$, then Assertion~(4)(b) is equivalent to $p\leq q$.}
\end{rem}

\begin{proof}[Proof of theorem~\ref{thm:mixed-exponents}](1) If $J_\varphi=\varnothing$, by continuity of $\varphi$ in a neighbourhood of $\overline{\D}$, there exist $R>1$ and $M<1$ such that
	$|\varphi(z)|\leq M$ for every $|z|\leq R$. Let $f=\sum_{k\geq 0}a_kz^k\in \ell_a^p$. Then
	\[
	|f\circ \varphi(z)| \leq \sum _{k\geq 0}|a_k| |\varphi(z)|^k \leq \frac{\Vert f \Vert _{\ell_a^p}}{1-M}.
	\]
	Writing $f\circ \varphi =\sum_{k\geq 0}b_k z^k$, we get, by Cauchy's inequalities,
	\[
	|b_k|\leq \frac{\Vert f \Vert _{\ell_a^p}}{1-M}R^{-k}.
	\]
	Thus, since $R>1$, we get, if $q\in[1,\infty)$,
	\[
	\Vert f\circ \varphi \Vert _{\ell_a^q} \leq \frac{\Vert f \Vert _{\ell_a^p}}{1-M}\left(\sum_{k\geq 0}R^{-qk}\right)^{1/q} \lesssim \Vert f \Vert _{\ell_a^p}.
	\]
	The case $q=\infty$ trivially follows.
	
	(2) Let $\varphi(z)=\mu z^n$ for some $|\mu|\leq 1$ and some $n\in \N$. We distinguish two cases. If $|\mu|<1$, then we are in the previous situation and $C_\varphi:\ell_a^p \to \ell_a^q$ is bounded for every $p,q \in [1,\infty]$. If $|\mu|=1$, then $\Vert C_\varphi(f)\Vert _{\ell_a^q}=\Vert f \Vert_{\ell_a^q}$ for every $f\in \ell_a^q$. Therefore, $C_\varphi: \ell_a^p \to \ell_a^q$ is bounded if and only if $\ell_a^p\subset \ell_a ^q$, which is equivalent to $p\leq q$.
	
 	(3) is a direct consequence of Corollary~\ref{cor-BF-intro} and the fact that $C_\varphi$ is bounded on $\ell_a^2$ for every $\varphi$.
 	
 	\medskip
	
	Let us now focus on Assertion~(4). Then $J_\varphi$ is finite and $r(\theta)\leq\sigma(\theta)<\infty$ for every  $\theta\in J_\varphi$. In particular, $D_\varphi\in \N$.
	
	We prove $\text{(b)}\Longrightarrow \text{(a)}$. Suppose first that $D_\varphi>0$, and put
	\[
	q_\varphi:=\frac{2(D_\varphi+1)}{D_\varphi+2}
	\quad\text{and}\quad
	p_\varphi:=\frac{2(D_\varphi+1)}{D_\varphi}.
	\]
	Thus $1<q_\varphi<2<p_\varphi<\infty$. Note that $p_\varphi=q_\varphi'$ is the conjugate exponent of $q_\varphi$. It is easily checked that for every $q\in[1,2]$, $\alpha_q\leq 0$ if and only if $q\geq q_\varphi$. In particular, $\alpha_{q_\varphi}=0$. Therefore, Corollary~\ref{main-thm-intro-1} gives the boundedness of the three following operators:
	\begin{equation}\label{eq:mixed-three-bounds}
		C_\varphi:\ell_a^1\to\ell_a^{q_\varphi},\quad C_\varphi:\ell_a^2\to\ell_a^2\quad \text{and}\quad C_\varphi:\ell_a^{p_\varphi}\to\ell_a^\infty.
	\end{equation}
	By Riesz--Thorin's theorem, interpolation between $[\ell_a^1,\ell_a^{q_\varphi}]$ and $[\ell_a^2,\ell_a^2]$ yields the boundedness of $C_\varphi:\ell_a^p \to \ell_a^q$, where, for $t\in [0,1]$,
	\[
	\frac1p=1-\frac t2 \quad \text{and}\quad \frac1q=\frac{1-t}{q_\varphi}+\frac t2=\frac{\frac1p+\frac{D_\varphi}{2}}{D_\varphi+1}=\frac{2+pD_\varphi}{2p(D_\varphi +1)}.
	\]
	Thus, using $\ell_a^{q_1}\subset \ell_a^{q_2}$ if $q_1\leq q_2$, for every $p\in [1,2]$, $C_\varphi:\ell_a^p \to \ell_a^q$ is bounded whenever
	\begin{equation}\label{cond-pq2}
	\frac{q}{2} \geq \frac{p(D_\varphi+1)}{pD_\varphi +2}.
	\end{equation}

	Similarly, interpolation between $[\ell_a^2,\ell_a^{2}]$ and $[\ell_a^{p_\varphi},\ell_a^{\infty}]$ in \eqref{eq:mixed-three-bounds} implies that, for every $p\in[2,p_\varphi]$, $C_\varphi:\ell_a^p\to \ell_a^q$ is bounded if
	\begin{equation}\label{cond-pq1}
	\frac{p}{2} \leq \frac{q(D_\varphi+1)}{qD_\varphi +2}.
	\end{equation}
	Note that conditions \eqref{cond-pq2} and \eqref{cond-pq1} are exactly those appearing in Assertion~(4) (with the usual convention when $p=\infty$ or $q=\infty$). Moreover, since $p_\varphi=\frac{2(D_\varphi+1)}{D_\varphi}$, a straightforward calculation shows that for every $p>p_\varphi$, there is no $q\in [1,\infty]$ for which \eqref{cond-pq1} holds. This proves (b)$\Longrightarrow$(a) in the case where $D_\varphi>0$.
	
	\medskip
	
	If $D_\varphi=0$, namely $r(\theta)=\sigma(\theta)$ for every $\theta\in J_\varphi$, then Corollary~\ref{cor-linfty-intro} yields the boundedness of $C_\varphi$ on $\ell_a^p$ for every $p\in[1,\infty]$. Since $\ell_a^p\subset\ell_a^q$ whenever $p\leq q$, we immediately obtain (b)$\Longrightarrow$(a).
	
	\medskip
	
	Let us now prove (a)$\Longrightarrow$(b). Assume first that $D_\varphi>0$ and choose $\theta_0\in J_\varphi$ such that $\sigma(\theta_0)-r(\theta_0)=D_\varphi$. Write $r=r(\theta_0)$ and $\sigma=\sigma(\theta_0)$. We keep the notation of	the proof of Theorem~\ref{thm:local-matrix-obstruction}, in particular
	\[
	T_N=N^{-1/\sigma},\quad	\Lambda_N=NT_N^r=N^{D_\varphi/\sigma},\quad	|I_N|\asymp|J_N|\asymp NT_N^{r-1}=T_N^{-1}\Lambda_N.
	\]
	Using the same $u^N$ as in that proof,	we have $\|u^N\|_{\ell^p}=|I_N|^{1/p}$ and $|(M_Nu^N)_k|\gtrsim\Lambda_N^{1/2}$ for all $k\in J_N'$, with $|J_N'|\asymp T_N^{-1}$. Thus, taking the $\ell^q$ norm instead of the	$\ell^p$ norm gives, for every $p,q\in[1,\infty]$,
	\begin{equation}\label{eq:mixed-first-lower}
		\begin{aligned}
			\|L_N\AA_\varphi P_N\|_{\ell^p\to\ell^q(\Z)}
			&\gtrsim\frac{T_N^{-1/q}\Lambda_N^{1/2}}
			{(T_N^{-1}\Lambda_N)^{1/p}}\\
			&=N^{\frac1\sigma\left(\frac{D_\varphi}{2}
				+\frac1q-\frac{D_\varphi+1}{p}\right)}.
		\end{aligned}
	\end{equation}
	Similarly, for the finite dimensional blocks $F_N$ and the same finite sequences $v^N$, one has $|J_N|^{1/q'}$ and	$|(F_N^*v^N)_n|\gtrsim\Lambda_N^{1/2}$ for every $n\in I_N'$, where
	$|I_N'|\asymp T_N^{-1}$. Thus
	\begin{equation}\label{eq:mixed-second-lower}
		\begin{aligned}
			\|L_N\AA_\varphi P_N\|_{\ell^p\to\ell^q(\Z)}
			&\geq\|F_N\|_{\ell^p(I_N)\to\ell^q(J_N)}\\
			&=\|F_N^*\|_{\ell^{q'}(J_N)\to\ell^{p'}(I_N)}\\
			&\gtrsim\frac{T_N^{-1/p'}\Lambda_N^{1/2}}
			{(T_N^{-1}\Lambda_N)^{1/q'}}\\
			&=N^{\frac1\sigma\left(\frac{D_\varphi+1}{q}
				-\frac1p-\frac{D_\varphi}{2}\right)}.
		\end{aligned}
	\end{equation}
	By Lemma~\ref{lem:uniform-localization}, the sequence $L_N$ is bounded on $\ell^q$ uniformly with respect to $N$. Therefore, if $C_\varphi:\ell_a^p\to\ell_a^q$ is bounded, then both exponents of $N$ in \eqref{eq:mixed-first-lower} and \eqref{eq:mixed-second-lower} are nonpositive, which is equivalent to the two inequalities in Assertion~(4). This proves (a)$\Longrightarrow$(b) when $D_\varphi>0$.
	
	\medskip
	
	Finally, it remains to consider the case $D_\varphi=0$. Fix $\theta_0\in J_\varphi$	and write, for $t$ close to $0$,
	\[
	\varphi(e^{i(\theta_0+t)})=\varphi(e^{i\theta_0})\exp(i\tau t+h(t)),
	\]
	with $\tau >0$ and $h(t)=O(t^2)$. Let $0<a<1<b$, $\eta\in C_c^\infty((a,b))$ nonnegative and equal to $1$ in a neighborhood of $1$, and set	$T_N=N^{-1}$ and $L_N=L_{T_N}$, as in Lemma~\ref{lem:uniform-localization}. For $\delta>0$ sufficiently small, set
	\[
	I_N=\{n\in\N:N\leq n\leq(1+\delta)N\},\quad\text{and}\quad J_N=\{k\in\N:|k-\tau N|\leq\delta N\},
	\]
	and consider the finitely supported sequence
	\[
	u^N(n)=
	\begin{cases}
		\varphi(e^{i\theta_0})^{-n},&n\in I_N,\\
		0,&n\notin I_N.
	\end{cases}
	\]
	Since $|\varphi(e^{i\theta_0})|=1$ and $|I_N|\asymp N$,	we have $\|u^N\|_{\ell^p}\asymp N^{1/p}$. Moreover, for every $k\in J_N$,
	\begin{equation}\label{eqfinalpq}
	e^{ik\theta_0}(L_N\AA_\varphi u^N)_k=\frac{1}{2\pi N}\int_a^b\eta(s)\sum_{n\in I_N}e^{nh(s/N)}e^{i(\tau n-k)s/N}ds.
	\end{equation}
	We have	$nh(s/N)=O(N^{-1})$ uniformly with respect to $n\in I_N$ and $s\in[a,b]$. Moreover,
	\[
	\left|\frac{(\tau n-k)s}{N}\right| 	\leq b(1+\tau)\delta,\quad n\in I_N,\quad k\in J_N.
	\]
	Upon decreasing $\delta$ if necessary, we may assume that $b(1+\tau)\delta\leq \pi/6$ and $\delta<\tau/2$. For all sufficiently large $N$,	the real part of \eqref{eqfinalpq} is then bounded from below by a positive constant, independent of $n\in I_N$, $k\in J_N$ and $s\in[a,b]$. Since $\eta\geq0$,
	\[
	|(L_N\AA_\varphi u^N)_k|\gtrsim\frac{|I_N|}{N}\gtrsim1,\quad k\in J_N.
	\]
	As $|J_N|\asymp N$, we conclude that
	\[
	\frac{\|L_N\AA_\varphi u^N\|_{\ell^q(\Z)}}
	{\|u^N\|_{\ell^p}}
	\gtrsim N^{1/q-1/p}.
	\]
	Since $L_N$ is uniformly bounded on $\ell^q$ by Lemma~\ref{lem:uniform-localization}, $C_\varphi:\ell_a^p \to \ell_a^q$ can be bounded only if $p\leq q$, which completes the proof.
\end{proof}

\begin{figure}[!htbp]
	\centering
	\begin{tikzpicture}[
		x=5.7cm,y=5.7cm,
		line cap=round,line join=round,
		axis/.style={->,line width=.55pt},
		boundary/.style={line width=1.05pt},
		guide/.style={densely dashed,line width=.45pt,black!48},
		tick/.style={line width=.5pt},
		critical/.style={circle,fill=black,inner sep=1.45pt},
		every node/.style={inner sep=1.7pt}
		]
		\pgfmathsetmacro{\Dvalue}{1}
		\pgfmathsetmacro{\pstar}{2*(\Dvalue+1)/\Dvalue}
		\pgfmathsetmacro{\qstar}{2*(\Dvalue+1)/(\Dvalue+2)}
		\pgfmathsetmacro{\xstar}{1/\pstar}
		\pgfmathsetmacro{\ystar}{1/\qstar}

		\fill[black!9] (\xstar,0)--(1,0)--(1,\ystar)--(.5,.5)--cycle;
		\draw[guide] (0,0)--(1,1);
		\node[font=\footnotesize,above left] at (.83,.83) {$p=q$};
		\draw[axis] (0,0)--(1.10,0) node[below right] {$1/p$};
		\draw[axis] (0,0)--(0,1.10) node[above left] {$1/q$};
		\draw[boundary] (\xstar,0)--(.5,.5)--(1,\ystar)--(1,0)--cycle;
		
		\foreach \v in {0,0.5,1}{
			\draw[tick] (\v,.006)--(\v,-.006);
			\draw[tick] (.006,\v)--(-.006,\v);
		}
		\draw[tick] (\xstar,.006)--(\xstar,-.006);
		\draw[tick] (.006,\ystar)--(-.006,\ystar);
		\node[font=\scriptsize,below left] at (0,0) {$0$};
		\node[font=\scriptsize,below] at (\xstar,0) {$1/p_\varphi$};
		\node[font=\scriptsize,below] at (.5,0) {$1/2$};
		\node[font=\scriptsize,below] at (1,0) {$1$};
		\node[font=\scriptsize,left] at (0,.5) {$1/2$};
		\node[font=\scriptsize,left] at (0,\ystar) {$1/q_\varphi$};
		\node[font=\scriptsize,left] at (0,1) {$1$};
		
		\node[critical] at (\xstar,0) {};
		\node[critical] at (.5,.5) {};
		\node[critical] at (1,\ystar) {};
		\node[font=\scriptsize,above left] at (.5,.5) {$(1/2,1/2)$};
		\node[font=\scriptsize,above] at (1,\ystar) {$(1,1/q_\varphi)$};
		\node[font=\footnotesize] at (.76,.19) {bounded};
	\end{tikzpicture}
	\caption{The shaded region, including its boundary, represents the set of those $(1/p,1/q)$ for which $C_\varphi:\ell_a^p\to\ell_a^q$ is bounded, when $D_\varphi>0$. Here $D_\varphi=1$, so that $q_\varphi=4/3$ and $p_\varphi=4$.
	\\If $D_\varphi = 0$, the corresponding region is the triangle $0\leq 1/q\leq 1/p \leq 1$.}
	\label{fig:pq-region-general}
\end{figure}
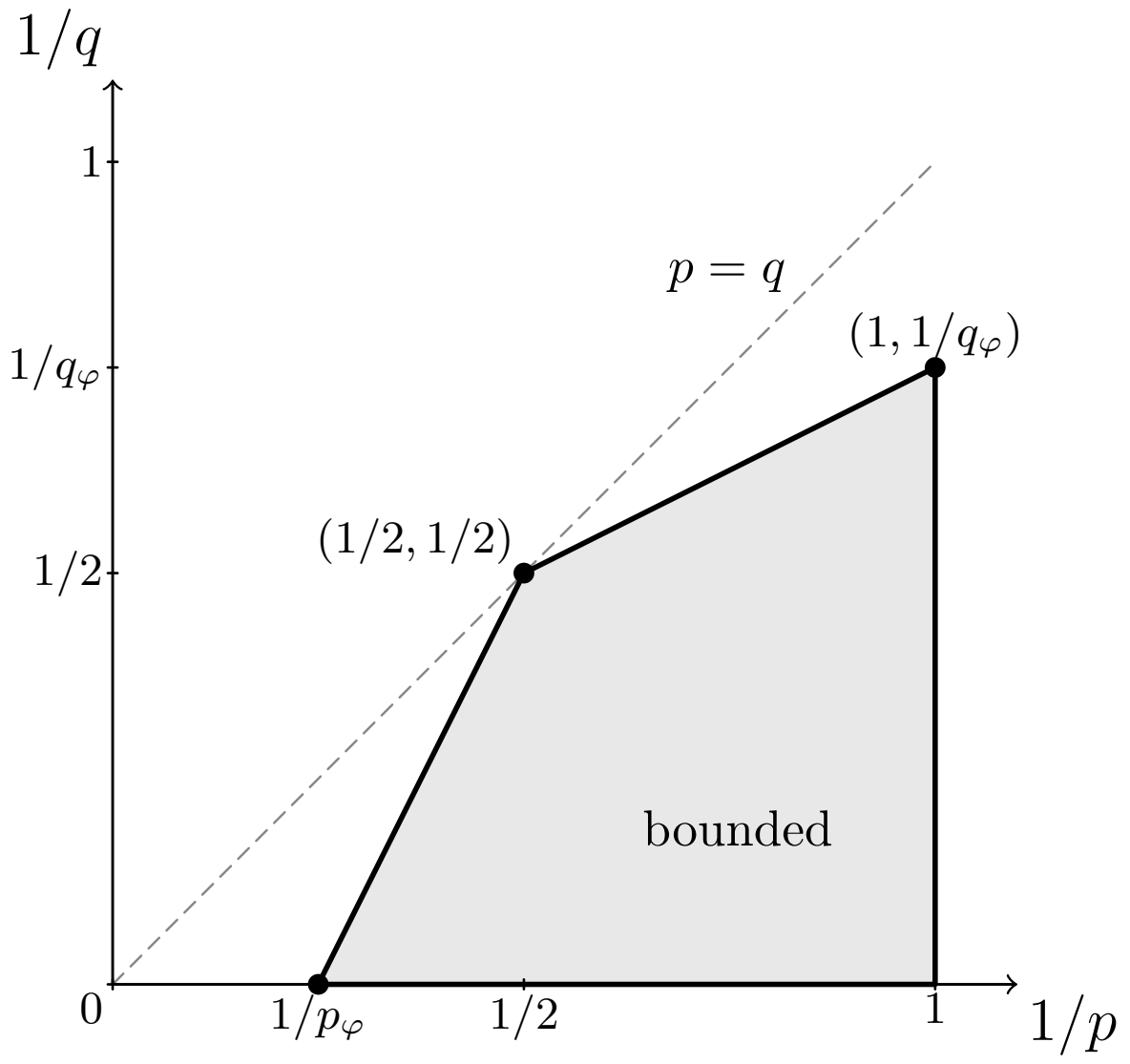

The boundedness region in Assertion~(4) is illustrated  in Figure~\ref{fig:pq-region-general}.

Let us give some examples.

\begin{examples}\label{examplesymbols}{\rm
		
		\noindent 1) Let $\varphi =\sum_{k\geq 0}a_k z^k$ be a symbol in $H(\overline{\D})$ such that $J_\varphi \neq \emptyset$. We assume that $a_k>0$ for every $k\in \N$. Observe that, by the triangle inequality, this implies $J_\varphi=\{0\}$. Then we claim that $r(0)=\sigma(0)$. Indeed, it is easily checked that, for $t$ close to $0$,
		\[
		\varphi(e^{it})=1+ic_1t-\frac{c_2}{2}t^2 + O(t^3),
		\]
		where
		\[
		c_1=\sum_{k\geq 0}ka_k\quad \text{and}\quad c_2=\sum_{k\geq 0}k^2a_k.
		\]
		Then, for $t$ close enough to $0$, we get
		\[
		\log(\varphi(e^{it}))=ic_1t-\frac{1}{2}(c_2-c_1^2)t^2 + O(t^3).
		\]
		The equality $\sum_{k\geq 0}a_k=1$ yields $c_2-c_1^2 = \sum_{k\geq0}a_k(k-c_1)^2$. Thus the assumption $a_k>0$ for at least two indices $k\geq 0$ implies $c_2-c_1^2 >0$, hence that $r(0)=\sigma(0)=2$.
		
		By Corollary~\ref{cor-linfty-intro}, the composition operator $C_\varphi$ is bounded on $\ell_a^p$ for every $p\in[1,\infty]$.
		
		\medskip
		
		\noindent 2) More generally, let $\varphi(z)=\sum_{k=0}^d a_kz^k$ be a nonconstant polynomial, with $d\geq1$ and $(a_1,\ldots,a_d)\ne(0,\ldots,0)$, and assume that
		\begin{equation}\label{cond-example-1}
			\sum_{k=0}^d|a_k|\leq1.
		\end{equation}
		Then $\varphi$ is a symbol and we claim that $C_\varphi$ is bounded on $\ell_a^p$ for every $p\in[1,\infty]$. Indeed, if the inequality is strict, then $J_\varphi=\emptyset$. If not, assume that $J_\varphi \neq \emptyset$ and, by Corollary~\ref{cor-linfty-intro}, it suffices to check that $D_\varphi = 0$. Fix $\theta \in J_\varphi$, and let $\zeta=e^{i\theta}$ and $\eta=\varphi(\zeta)$. The equality in the triangle inequality implies
		$\overline\eta a_k\zeta^k=|a_k|$ for every $k$, hence
		\[
		\overline\eta\,\varphi(\zeta e^{it})=\sum_{k=0}^d|a_k|e^{ikt}=1+i\mu t-\frac{\nu}{2}t^2+O(t^3),
		\]
		where, as in the first example,
		\[
		\mu=\sum_{k=0}^d k|a_k|\quad\text{and} \quad	\nu=\sum_{k=0}^d k^2|a_k|.
		\]
		Therefore,
		\[
		\log\!\bigl(\overline\eta\,\varphi(\zeta e^{it})\bigr)
		=i\mu t-\frac{\nu-\mu^2}{2}t^2+O(t^3).
		\]
		Now, one can easily check (again as in 1)) that, since $\sum_k|a_k|=1$,
		\begin{equation*}
			\nu-\mu^2 = \sum_{0\leq k<j\leq d}|a_k||a_j|(k-j)^2>0,
		\end{equation*}
		if two coefficients are nonzero. Consequently, $r(\theta)=\sigma(\theta)=2$ at every contact point. If only one coefficient is nonzero, then $\varphi$ is monomial
		
		\medskip
		
		In particular, this shows that every nonconstant affine symbol induces a bounded composition operator on $\ell_a^p$ for every $p\in[1,\infty]$.
		
		\medskip
		
		\noindent 3) Let $\varphi=\sum_{k=0}^d a_kz^k$ be such that	$\varphi(1)=1$ and $J_\varphi=\{0\}$. In this case,	$\varphi'(1)=\sum_{k=1}^d ka_k>0$ is real. Setting $M_2=\sum_{k=1}^d k^2a_k$, we may check that
		\[
		\log\varphi(e^{it})	=i\varphi'(1) t-\frac{M_2-\varphi'(1)^2}{2}t^2+O(t^3).
		\]
		Therefore,
		\begin{equation}\label{example-cond}
			\Re\!\left(\sum_{k=1}^d k^2a_k\right)>\varphi'(1)^2
		\end{equation}
		implies $r(0)=\sigma(0)=2$, which quarantees the boundedness of $C_\varphi$ on every $\ell_a^p$ for every $p\in[1,\infty]$, by Corollary~\ref{cor-linfty-intro}.
		
		\medskip
		
		We point out that for $d=2$ and $d=3$, condition~\eqref{example-cond} is equivalent to
		\[
		\Re(a_1+4a_2)>(a_1+2a_2)^2,\quad \text{and} \quad\Re(a_1+4a_2+9a_3)>(a_1+2a_2+3a_3)^2,
		\]
		respectively (the quantities on the right are equal to $\varphi'(1)^2$ and thus are real).
		
		\medskip
		
		\noindent 4) Let us consider the polynomials
		\[
		\varphi_2(z)=\frac{3+6z-z^2}{8}\quad \text{and} \quad \varphi_3(z)=\frac{-1+9z+9z^2-z^3}{16}.
		\]
		Straightforward computations show that $J_{\varphi_2}=J_{\varphi_3}=\{0\}$. Moreover,
		\[
		\sum_k|\widehat{\varphi_2}(k)| =  \sum_k|\widehat{\varphi_3}(k)| = \frac54.
		\]
		Further,
		\[
		\varphi_2'(1)=\frac12,\quad	\sum_k k^2\widehat{\varphi_2}(k)=\frac14,	\quad\varphi_3'(1)=\frac32\quad \text{and}\quad \sum_k k^2\widehat{\varphi_3}(k)=\frac94.
		\]
		Thus, neither $\varphi_2$ nor $\varphi_3$ satisfy the assumptions of the two previous examples. Now, a calculation, which is left to the reader, shows that
		\begin{align*}
			\log\varphi_2(e^{it}) & =\frac{i}{2}t+\frac{i}{16}t^3-\frac{3}{128}t^4+O(t^5),\\
			\log\varphi_3(e^{it})&=\frac{3i}{2}t-\frac{3}{128}t^4+O(t^5).
		\end{align*}
		Hence $r_{\varphi_2}(0)=3<\sigma_{\varphi_2}(0)=4$, whereas	$r_{\varphi_3}(0)=\sigma_{\varphi_3}(0)=4$. In particular, $C_{\varphi_2}$ is bounded on $\ell_a^p$ if and only if $p=2$ by Theorem~\ref{main-thm-rlesssigma}, while $C_{\varphi_3}$ is bounded on $\ell_a^p$ for every	$1\leq p\leq\infty$, by Corollary~\ref{cor-linfty-intro}.
		
		The last assertion also shows that neither condition~\ref{cond-example-1} nor condition~\ref{example-cond} is necessary for the boundedness of $C_\varphi$ on $\ell_a^p$, $p\neq 2$.
	}
\end{examples}

\section{Extensions to weighted Beurling--Sobolev spaces $\ell_a^p(w)$}\label{section-1-Charac}

In this section, we discuss to which extent the results of the previous sections extend to the weighted setting, that is to $\ell_a^p(w)$ spaces. For convenience, we focus on composition operators $C_\varphi$ acting on $\ell_a^p(w)$ (and not between two different weighted Beurling--Sobolev spaces).

In \cite{LLQR-l2weight-2025,LLQR-lpa-2025}, as recalled in the introduction, it appears that slow oscillation of the weight $w$ may play an important role in the Hilbertian setting: the composition operator induced by a single non-trivial automorphism is bounded on $\ell_a^2(w)$ if and only if $w$ is slowly oscillating; on the contrary, in the non-Hilbertian setting, no such composition operator is bounded on $\ell_a^p(w)$, whatever $w$ and $p\neq 2$. We shall see that, in general, things are more complicated. A thorough study of composition acting on $\ell_a^p(w)$ for general weights can be found in \cite{CharpentierEspoullierZarouf2026-1}.

\medskip

First of all, let us treat the simple case $J_\varphi = \emptyset$. Recall that the latter holds if and only if $\varphi(\overline{\D})\subset \D$. In this case, the composition operator $C_{\varphi}$ acts regularly on $\ell_a^p(w)$ for the whole range of $p$, for every \textit{reasonable} weight:

\begin{prop}\label{prop-jphiempty}Let $w$ be a weight such that $w_k^{1/k}\to 1$ as $k\to \infty$. For every symbol $\varphi \in H(\overline{\D})$, if $J_{\varphi}=\emptyset$, then $C_{\varphi}$ maps $\ell_a^p(w)$ into itself, for every $p \in [1,\infty]$.
\end{prop}

\begin{proof}By the Riesz--Thorin interpolation theorem, it is enough to consider the cases $p=1$ and $p=\infty$. By Proposition~\ref{prop-key-boundedness-carac}, it suffices to check that
	\begin{equation}\label{eq-Jphiempty}
		\sup_{n\geq 0} \sum_{k\geq 0} |\widehat{\varphi^n}(k)|\frac{w_k}{w_n} < \infty \quad \text{and} \quad \sup_{k\geq 0} \sum_{n\geq 0} |\widehat{\varphi^n}(k)|\frac{w_k}{w_n} < \infty.
	\end{equation}
	By the assumption, there exist $M\in [0,1)$ and $R>1$ such that $|\varphi(z)|\leq M$ for every $z$ with $|z|\leq R$. By Cauchy's estimates, it follows that, for every $n,k\in \N$,
	\[
	|\widehat{\varphi^n}(k)|\leq \frac{M^n}{R^k}.
	\]
	Choose $a,b$ with $M<a<1<b<R$. Since $w_j^{1/j}\to1$, there exists $C\geq1$ such that
	\[
	C^{-1}a^j\leq w_j\leq Cb^j,\qquad j\geq0.
	\]
	Together with Cauchy's estimate, this gives two convergent geometric series, and hence \eqref{eq-Jphiempty}.
\end{proof}

Let us now focus on slowly oscillating weights.

\subsection{Slowly oscillating weights}

We recall the definition.
\begin{defi}A weight $w$ is said to be slowly oscillating if, for every $0<J<1<L<\infty$, there exists $c> 0$ such that, for every $n\in \N$ and every $Jn\leq k \leq Ln$, one has
	\[
	c^{-1}w_{n}\leq w_k\leq cw_{n}.
	\]
\end{defi}

This definition coincides with Definition~2.2 in \cite{LLQR-l2weight-2025}. Typical examples of slowly oscillating weights are
\[
w_k=(k+1)^{\alpha}\bigl(\log(e+k)\bigr)^{\beta},
\]
where $\alpha,\beta\in\mathbb R$. We refer to \cite{LLQR-l2weight-2025} for further examples.

By \cite[Proposition~2.5]{LLQR-l2weight-2025}, every slowly oscillating weight $w$ admits polynomial upper and lower bounds: there exist $\alpha,\beta>0$ such that
\begin{equation}\label{poly-esti-weight}
	k^{-\alpha}\lesssim w_k\lesssim k^{\beta}.
\end{equation}
Moreover, one can easily check that, if $w$ is a slowly oscillating weight, then $w_k^{1/k}  \to 1$ as $k\to \infty$. Therefore, by Proposition~\ref{prop-jphiempty}, if $\varphi \in H(\overline{\D})$ satisfies $J_\varphi = \emptyset$, then $C_\varphi$ is bounded on $\ell_a^p(w)$ for every $p\in[1,\infty]$ and every slowly oscillating weight $w$. Note that the weight defined by $w_k=e^{-(\log(k+1))^2}$ satisfies $w_k^{1/k}  \to 1$, but it is not slowly oscillating.

\medskip

The main result of \cite{LLQR-l2weight-2025} implies that, if the weight  $w$ is slowly oscillating, then the boundedness of $C_\varphi$ on $\ell_a^2(w)$ is equivalent to that on $\ell_a^2$, for every symbol $\varphi \in H(\overline{\D})$. More precisely, \cite[Theorems~4.6 and 4.9]{LLQR-l2weight-2025} read as follows.

\begin{thm}[\cite{LLQR-l2weight-2025}]\label{thm:Hbar-SO}Let $w$ be a weight. Then the following assertions are equivalent.
	\begin{enumerate}
		\item $w$ is slowly oscillating;
		\item $C_\varphi$ is bounded on $\ell_a^2(w)$ for every $\varphi\in H(\overline{\D})$;
		\item $C_\varphi$ is bounded on $\ell_a^2(w)$ for some automorphism $\varphi$ of $\D$, with $\varphi(0)\neq 0$.
	\end{enumerate}
\end{thm}

The next results show that, for a symbol in $H(\overline{\D})$, the weighted setting exhibits, also for $p\neq 2$, the same phenomenon as the unweighted setting, whenever the weight is slowly oscillating.

\begin{thm}\label{prop-key-boundedness-carac-2}Let $w$ be a slowly oscillating weight and let $\varphi \in H(\overline{\D})$ be a symbol. For every $p\in[1,\infty)$,
	\[
		\sup_{n\geq 0} \sum_{k\geq 0} |\widehat{\varphi^n}(k)|^p < \infty \quad \Longleftrightarrow \quad \sup_{n\geq 0} \sum_{k\geq 0}|\widehat{\varphi^n}(k)|^p\frac{w_k^p}{w_n^p} <\infty
	\]
	and, for every $p\in[1,\infty)$,
	\[
		\sup_{k\geq 0} \sum_{n\geq 0} |\widehat{\varphi^n}(k)|^{p} < \infty \quad \Longleftrightarrow \quad \sup_{k\geq 0} \sum_{n\geq 0}|\widehat{\varphi^n}(k)|^p\frac{w_k^p}{w_n^p} <\infty.
	\]
	In particular, Theorem~\ref{thm-intro-main} and Corollaries~\ref{main-thm-intro-1}, \ref{cor-BF-intro} and \ref{cor-linfty-intro} remain true if $\ell_a^p$ is replaced with $\ell_a^p(w)$.
\end{thm}

\begin{proof}We shall only prove the second equivalence. The first one is proved the same way. \begin{blueadd}
	For $k\geq0$, set
	\[
	S(k):=\sum_{n\geq0}|\widehat{\varphi^n}(k)|^p,
	\qquad
	S_w(k):=\sum_{n\geq0}|\widehat{\varphi^n}(k)|^p\frac{w_k^p}{w_n^p}.
	\]
	By Lemma~\ref{lem-elementary}, the terms with $n<Jk$ or $n>Lk$ in $S(k)$ are exponentially small. On the other hand, slow oscillation gives $w_n\asymp w_k$ for $Jk\leq n\leq Lk$. Hence, for some $\rho\in(0,1)$,
	\[
	S(k)\lesssim S_w(k)+\rho^k,
	\]
	which proves that $\sup_k S_w(k)<\infty$ implies $\sup_k S(k)<\infty$.

	Conversely, the polynomial bounds on $w$, together with Lemma~\ref{lem-elementary} using suitable powers of $n$, show that the weighted tails $n<Jk$ and $n>Lk$ are exponentially small. The same comparison $w_n\asymp w_k$ on $Jk\leq n\leq Lk$ therefore gives
	\[
	S_w(k)\lesssim S(k)+\rho^k.
	\]
	Thus $\sup_k S(k)<\infty$ implies $\sup_k S_w(k)<\infty$.
\end{blueadd}
	This yields the second equivalence.
\end{proof}

The following slowly oscillating weights version of Theorem~\ref{thm:local-matrix-obstruction} holds also true.

\begin{thm}\label{thm:weighted-quant}
	Let $w$ be a slowly oscillating weight, let $p\in[1,\infty]$, and let $\varphi\in H(\overline{\D})$ be a symbol. Assume that there exists $\theta_0\in J_\varphi$ such that
	\[
	r:=r(\theta_0)<\sigma:=\sigma(\theta_0)<\infty.
	\]
	Let $P_N$ be the coordinate projections from Theorem~\ref{thm:local-matrix-obstruction}. Then there exist constants $0<c<C<\infty$ and, for every sufficiently large $N$, a coordinate projection $R_N$ of $\ell^p$ onto $[cN,CN] \cap \N$, such that
	\[
	\norm{R_N\AA_\varphi  P_N}_{\ell^p(w)\to\ell^p(w)} \gtrsim N^{\frac{\sigma-r}{\sigma}\left|\frac1p-\frac12\right|}.
	\]
	In particular, $C_\varphi$ is unbounded on $\ell_a^p(w)$ if $p\neq 2$.
\end{thm}

\begin{proof}We have to distinguish the cases $p<\infty$ and $p=\infty$. For convenience, we set $\tilde{p}=p$ if $p<\infty$ and $\tilde{p}=1$ if $p=\infty$. By Lemma~\ref{lem-elementary}, there exist constants $0<J<1<L<\infty$
	and $s\in(0,1)$ such that, for every sufficiently large $n$,
	\begin{equation}\label{4}
		\sum_{k<Jn}|\widehat{\varphi^n}(k)|^{\tilde{p}} +\sum_{k>Ln}|\widehat{\varphi^n}(k)|^{\tilde{p}} \leq s^n.
	\end{equation}
	In the sequel, the notation of the proof of Theorem~\ref{thm:local-matrix-obstruction} is retained. In particular
	\[
	I_N=\{n\in\N:|n-N|\le \delta NT_N^{r-1}\},
	\]
	with $T_N=N^{-1/\sigma}$. The matrix $P_N$ is the coordinate projection of $\ell^p$ onto $I_N$, and $L_N=L_{T_N}$ is the convolution operator defined in \eqref{defi-conv-op-L_T}. Recall that $P_N$ and $L_N$ are uniformly bounded in norm.
	By Theorem~\ref{thm:local-matrix-obstruction},
	\[
	\Vert L_N\AA_\varphi P_N\Vert_{\ell^p\to\ell^p(\Z)}\gtrsim \Lambda_N^{\left|\frac1p-\frac12\right|}.
	\]
	Upon increasing $N$ if necessary, we may and shall assume
	\begin{equation}\label{5}
		I_N\subset [N/2,2N].
	\end{equation}
	Set
	\begin{equation*}\label{6}
		\mathcal K_N:=\bigl\{k\in\N:JN/2\le k\le 2LN\bigr\},
	\end{equation*}
	and let $R_N$ be the coordinate projection of $\ell^p$ onto $\mathcal K_N$. We shall first prove that, for all $N$ large enough,
	\begin{equation}\label{12}
		\Vert L_NR_N\AA_\varphi P_N\Vert_{\ell^p\to \ell^p(\Z)}
		\gtrsim \Lambda_N^{\left|\frac1p-\frac12\right|}.
	\end{equation}
	By the triangle inequality,
	\[
	\begin{aligned}
		\Vert L_NR_N\AA_\varphi P_N\Vert_{\ell^p\to \ell^p(\Z)} &\geq \Vert L_N\AA_\varphi P_N\Vert_{\ell^p\to \ell^p(\Z)}\\
		&\quad- \Vert L_N(I-R_N)\AA_\varphi P_N\Vert_{\ell^p\to \ell^p(\Z)}.
	\end{aligned}
	\]
	Let us show that $\Vert (I-R_N)\AA_\varphi P_N\Vert_{\ell^p\to \ell^p}$ is \textit{small} for every $N$ large enough.
	Note that, if $n\in I_N$ and $k\notin\mathcal K_N$, then \eqref{5} implies either $k<Jn$ or $k>Ln$. Therefore \eqref{4} gives, uniformly for $n\in I_N$,
	\begin{equation}\label{7}
		\sum_{k\notin\mathcal K_N}|\widehat{\varphi^n}(k)|^{\tilde{p}}\le s^n\le s^{N/2}.
	\end{equation}
	Assume first that $p<\infty$. For $a=(a_n)\in\ell^p(\N)$, \eqref{7} thus yields
	\begin{align*}
		\Vert(I-R_N)\AA_\varphi P_N(a) \Vert_{\ell^p}
		&\leq \sum_{n\in I_N}|a_n|
		\left(\sum_{k\notin\mathcal K_N}|\widehat{\varphi^n}(k)|^p\right)^{1/p}\\
		&\leq s^{N/(2p)}\sum_{n\in I_N}|a_n|\\
		&\leq |I_N|^{1/p'}s^{N/(2p)}\Vert a\Vert_{\ell^p},
	\end{align*}
	where $p'$ is the conjugate exponent of $p$, with the convention $|I_N|^{1/p'}=1$ if $p=1$. Since $|I_N|\lesssim NT_N^{r-1}= o(N)$, one has $|I_N|^{1/p'}\lesssim N$, hence
	\[
	N s^{N/(2p)}\lesssim s^{N/(2p+1)}.
	\]
	It follows
	\begin{equation}\label{9}
		\Vert(I-R_N)\AA_\varphi P_N\Vert_{\ell^p\to\ell^p}\lesssim s^{N/(2p+1)}.
	\end{equation}
	If $p=\infty$, then, for $a\in\ell^\infty(\N)$,
	\begin{align*}
		\Vert(I-R_N)\AA_\varphi P_N(a)\Vert_{\ell^\infty}
		&\leq \sum_{n\in I_N}|a_n|
		\sup_{k\notin\mathcal K_N}|\widehat{\varphi^n}(k)|\\
		&\leq \sum_{n\in I_N}|a_n|
		\sum_{k\notin\mathcal K_N}|\widehat{\varphi^n}(k)|\\
		&\le |I_N|s^{N/2}\Vert a \Vert_{\ell^\infty},
	\end{align*}
	and, since $|I_N| = o(N)$, we similarly obtain
	\[
	\Vert(I-R_N)\AA_\varphi P_N\Vert_{\ell^\infty\to\ell^\infty}\lesssim N s^{N/2}\lesssim s^{N/3}.
	\]
	The preceding estimates and the fact that $L_N$ is uniformly bounded show that
	$\Vert L_N(I-R_N)\AA_\varphi P_N \Vert$ is exponentially small in $N$. By Theorem~\ref{thm:local-matrix-obstruction}, the definition of $\Lambda_N$ yield \eqref{12}. Using again the uniform boundedness of $L_N$, we get
	\begin{equation}\label{13}
		\Vert R_N\AA_\varphi P_N\Vert_{\ell^p\to\ell^p}\gtrsim \Lambda_N^{\left|\frac1p-\frac12\right|},
	\end{equation}
	for every $N$ large enough.
	We now turn to the weighted norm. By the definition of $\mathcal K_N$, and since $w$ is slowly oscillating, there exists $C_w\geq1$, independent of $N$, such that
	\[
	C_w^{-1}w_N\le w_j\le C_ww_N,
	\quad j\in \mathcal K_N.
	\]
	Consequently, for every sequence $a$ supported in $\mathcal{K}_N$,
	\[
	C_w^{-1}w_N\Vert a\Vert_{\ell^p} \leq \Vert a\Vert_{\ell^p(w)}\leq C_ww_N\Vert a \Vert_{\ell^p}.
	\]
	Since $I_N\subset\mathcal K_N$ by \eqref{5} and $R_N\AA_\varphi P_N(a)$ is supported in $\mathcal K_N$ for every $a\in \ell^p$, \eqref{13}  implies
	\[
	\Vert R_N\AA_\varphi P_N\Vert_{\ell^p(w)\to\ell^p(w)}\geq C_w^{-2}\Vert R_N\AA_\varphi P_N\Vert_{\ell^p\to\ell^p}\gtrsim \Lambda_N^{\left|\frac1p-\frac12\right|},
	\]
	which proves the asserted estimate.
	Since $P_N$ and $R_N$ are uniformly bounded, the last assertion of the theorem directly follows because $\Lambda_N\to\infty$.
\end{proof}

The following theorem finally asserts that slowly oscillating weights do not influence the boundedness of composition operators on $\ell_a^p(w)$ spaces, compared to the unweighted setting.

\begin{thm}\label{thm:qualitative}Let $w$ be a slowly oscillating weight and let $\varphi\in H(\overline{\D})$ be a symbol. Then, for every $p\in[1,\infty]$,
	\[
		C_\varphi:\la^p(w)\to\la^p(w)\text{ bounded}\quad\Longleftrightarrow\quad C_\varphi:\la^p\to\la^p\text{ bounded}.
	\]
\end{thm}

\begin{proof}Let $w$ be slowly oscillating. For $p=2$, the asserted equivalence follows from Theorem~\ref{thm:Hbar-SO}. In the rest of the proof, we assume that $p\neq 2$. 
	
	If $J_\varphi = \emptyset$, then Proposition~\ref{prop-jphiempty} and the fact that $w$ slowly oscillating satisfies $w_k^{1/k}\to 1$ as $k\to \infty$, imply that $C_\varphi$ is bounded on every $\ell_a^p(w)$, $p\in[1,\infty]$. The same conclusion holds if $\varphi$ is a monomial.
	
	Suppose that $\varphi$ is a non-monomial finite Blaschke product. Then, by Corollary~\ref{cor-BF-intro}, Proposition~\ref{prop-key-boundedness-carac} and Theorem~\ref{prop-key-boundedness-carac-2}, $C_B$ is not bounded from $\ell_a ^1(w)$ to $\ell_a^p(w)$, $p\in [1,2)$, nor from $\ell_a^p(w)$ to $\ell_a^\infty(w)$, $p\in(2,\infty]$. Therefore,
	\[
	{C_B}:\, \ell_a^p(w)\to \ell_a^p(w)
	\]
	is not bounded, hence the asserted equivalence.
	
	We now consider the case where $\varphi$ is not monomial nor a finite Blaschke product, and where $J_\varphi \neq \emptyset$ and $r(\theta) = \sigma(\theta)$ for every $\theta \in J_\varphi$. By Corollary~\ref{cor-linfty-intro}, $C_\varphi$ is bounded on $\ell_a^1$ and on $\ell_a^\infty$. Thus, by Proposition~\ref{prop-key-boundedness-carac} and Theorem~\ref{prop-key-boundedness-carac-2}, the same holds true on $\ell_a^1(w)$ and $\ell_a^\infty(w)$. The Riesz-Thorin interpolation theorem (valid for weighted $\ell_a^p(w)$ spaces) implies that $C_\varphi$ is bounded on $\ell_a^p(w)$.
	
	The remaining case, corresponding to $\varphi$ not monomial, not a finite Blaschke product, $J_\varphi \neq \emptyset$ and $r(\theta) < \sigma(\theta)$ for some $\theta \in J_\varphi$, follows from Theorem~\ref{thm:weighted-quant}
\end{proof}

\subsection{Non slowly oscillating weights}\label{section-2-non-slowly-oscillation}
We shall see that none of the main results of the previous sections hold for general weights.

\subsubsection{A weight $w$ and a symbol $\varphi$ such that $r(\theta)=\sigma(\theta)=2$ for every $\theta \in J_{\varphi}$, $C_\varphi$ does not map $\ell_a^1(w)$ into itself, but $C_\varphi$ maps $\ell_a^{\infty}(w)$ into itself}\label{sub-counterexample-1}

Let us consider the symbol
\[
\varphi_0(z)=\frac{z+z^2}{2}.
\]

For every $n\in \N$, one has $\varphi_0^n(z)=2^{-n}z^n(1+z)^n$, so that
\begin{equation*}\label{eq:binomial-matrix}
	\widehat{\varphi_0^n}(k)=
	\begin{cases}
		2^{-n}\binom{n}{k-n}& \text{if }n\le k\le2n\\
		0&\text{otherwise}.
	\end{cases}
\end{equation*}
Moreover, one easily checks that $J_{\varphi_0}=\{0\}$ and $\log(\varphi_0(e^{i t}))=\frac{3i t}{2}- \frac{t^2}{8} + o(t^2)$ as $t\to 0$. Thus, $r(0)=\sigma(0)=2$ and $\alpha_1(\varphi_0)=0$.

Now, set $w =(w_k)_k$ with
\[
w_k=
\begin{cases}
	3^{-\frac{j}{2}} & \text{if } k=3^j\text{ for some }j\\
	1 &\text{otherwise}.
\end{cases}
\]
One has $\liminf w_k^{1/k}=1$, hence $w$ is a weight. It is not slowly oscillating, since
\[
\frac{w_{3^j+1}}{w_{3^j} }= \sqrt{3^j} \to \infty.
\]

The following proposition shows that none of the equivalences in Corollaries~\ref{main-thm-intro-1}--\ref{cor-linfty-intro} hold for general weights.

\begin{prop}With the above notation, we have
	\[
	J_{\varphi_0}=\{0\},\quad r(0)=\sigma(0)=2 \quad \text{and} \quad \alpha_1(\varphi_0)=0,
	\]
	However
	\[
	\sup_{k\geq0}\sum_{n\geq 0}\widehat{\varphi_0^n}(k) \frac{w_k}{w_n} <\infty \quad \text{and}\quad \sup_{n\geq0}\sum_{k\geq 0}\widehat{\varphi_0^n}(k) \frac{w_k}{w_n} =\infty.
	\]
	Thus, $C_{\varphi_0}$ is bounded on $\ell_a^{\infty}(w)$, but it is not bounded on $\ell_a^1(w)$.
\end{prop}

\begin{proof}Observe that, for $z\in \D$, one has
	\[
	\sum_{k\geq 0}\left(\sum_{n\geq 0}\widehat{\varphi_0^n}(k)\right) z^k  = \sum_{n\geq 0}\varphi_0^n(z)= \frac{1}{1-\varphi_0(z)}=\frac{2}{3}\sum_{k\geq 0}z^k + \frac{1}{3}\sum_{k\geq 0}(-1/2)^kz^k.
	\]
	Thus
	\[
	\sum_{n\geq 0}\widehat{\varphi_0^n}(k)={\frac{2}{3}+\frac{1}{3}(-\frac12)^k\leq 1}.
	\]	
	Fix $k\geq0$. For each interval $[k/2,k]$, there exists at most one $j\geq 0$ such that $k/2\leq 3^j \leq k$. If such a $j$ exists, then
	\[
	\sum_{n\geq 0} \widehat{\varphi_0^n}(k) \frac{w_k}{w_n} \leq 1 + 3^{j/2} \widehat{\varphi_0^{3^j}}(k).
	\]
	Using the standard estimate
	\[
	\max_{0\leq m\leq n}2^{-n}\binom{n}{m} \leq \frac{C}{\sqrt{n}}
	\]
	for some constant $C>0$, we get $\widehat{\varphi_0^{3^j}}(k) \leq C3^{-j/2}$ and then
	\[
	\sum_{n\geq 0} \widehat{\varphi_0^n}(k) \frac{w_k}{w_n} \leq 1 + C.
	\]
	If there is no such $j$, then the latter quantity is at most $1$.
	
	For the other quantity, consider the interval $[3^j,2\cdot3^j]$. This interval contains no other $3^i$, $i\neq j$. Since $\sum_{k\geq 0}\widehat{\varphi_0^n}(k)=\varphi_0(1)^n=1$ and $\widehat{\varphi_0^n}(n)=2^{-n}$, we have
	\[
	\sum_{k\geq 0}\widehat{{\varphi_0}^{3^j}}(k) \frac{w_k}{w_{3^{j}}} = 3^{j/2}(1-2^{-3^j}) + 2^{-3^j} \to \infty,
	\]
	as $j\to \infty$.
	
	The last assertion follows from Proposition~\ref{prop-key-boundedness-carac}.
	
\end{proof}

\begin{rem}{\rm Setting $w^{-1}=(w_k^{-1})_{k\geq 0}$, one can similarly show that $C_{\varphi_0}$ is bounded on $\ell_a^1(w^{-1})$ while it is not on $\ell_a^{\infty}(w^{-1})$.
	}
\end{rem}

\subsubsection{A class of non-slowly oscillating weights and finite Blaschke product composition operators bounded on every $\ell^p(w)$}
For $\lambda \in \D$, recall that $b_{\lambda}$ denotes the single Blaschke factor associated with $\lambda$. According to the result of Lef\`evre \emph{et al.} \cite{LLQR-lpa-2025}, for every $\lambda\in\mathbb D\setminus\{0\}$, the composition operator $C_{b_\lambda}$ is unbounded on $\ell_a^p(w)$ whenever $p\neq 2$, regardless of the choice of the weight $w$. By  \cite{LLQR-l2weight-2025}, $C_{b_\lambda}$ is bounded on $\ell_a^2(w)$ for every $\lambda \in \D$ if and only if $w$ is slowly oscillating (see Theorem~\ref{thm:Hbar-SO}). This may suggest that the slowly oscillating assumption in Corollary~\ref{cor-BF-intro} is merely technical and, for instance, that the conclusion of \cite{LLQR-lpa-2025} should continue to hold for arbitrary weights and all non-monomial finite Blaschke products. Somewhat counterintuitively, this is not the case, as the next result shows.

\medskip

For a finite Blaschke product $B$ and an integer $d\geq1$, we write $B_d$ for the finite Blaschke product defined by
\[
B_d(z)=B(z^d),\quad z\in\mathbb D.
\]

\begin{thm}\label{bounded-CB-not-SO}
	Let $w=(w_n)_{n\geq 0}$ be a weight. Assume that $w$ is decreasing and that there exists an integer $\mu \geq 0$ such that
	\begin{equation}\label{assumption-weight-bddness}
		\sup_{n\geq 0}n\frac{w_{\mu n}}{w_n} <\infty.
	\end{equation}
	Then, for every $p\in [1,\infty]$ and every finite Blaschke product $B$,
	\begin{enumerate}
		\item\label{tiem1thmlast} there exists an integer $d \geq 1$ such that the composition operator $C_{B_d}$  is bounded on $\ell_a^p(w)$.
		\item\label{tiem2thmlast} if $B(0)=0$, then there exists an integer $d \geq 1$ such that the composition operator $C_{B^d}$ is bounded on $\ell_a^p(w)$.
	\end{enumerate}
\end{thm}

\begin{rem}\label{remark-non-so-weight}{\rm Note that \eqref{assumption-weight-bddness} is equivalent to the condition: for any $\beta>0$, there exists $\mu \in \N$ such that
		\[
		\sup_{n \geq 0}n^{\beta}\frac{w_{\mu n}}{w_n} < \infty.
		\]
		Let $\gamma >0$. Then, the weight $w_n=\exp(-(\log(n+1))^{\gamma})$, $n\in \N$, satisfies \eqref{assumption-weight-bddness} if and only if $\gamma\geq 2$.
	}
\end{rem}

\begin{proof}[Proof of Theorem \ref{bounded-CB-not-SO}]Let us prove \eqref{tiem1thmlast}. Fix $w$ and $\mu$ as in the statement. Let $B$ be a finite Blaschke product. By Lemma \ref{lem-elementary}, there exist $J\in (0,1)$, $q\in (0,1)$ and $C_0\geq1$ such that for every $n\geq0$,
	\begin{equation}\label{eq-reuse-lemma}
		\sum_{k=0}^{Jn}|\widehat{B^n}(k)| \leq C_0q^n.
	\end{equation}
	Let $n_0\in \N$ be such that $\lfloor Jn\rfloor\neq0$ for every $n\geq n_0$, and let $d\in \N$ be such that, for every $n\geq n_0$, $d\lfloor Jn\rfloor \geq \mu n$, and set $B_d (z)= B(z^d)$. By the Riesz-Thorin interpolation theorem, it suffices to check that $C_{B_d}$ acts boundedly on $\ell_a^p(w)$ for $p=1$ and $p=\infty$.
	
	Let us start with $p=1$. 
	Observe that $C_{B_d}$ is bounded on $\ell_a^1(w)$ if and only if
	\[
	\sup_{n\geq n_0} \frac{\sum _{k\geq 0}|\widehat{B^n}(k)|w_{dk}}{w_n}<\infty.
	\]
	Recall that $w$ is decreasing. Thus, by \eqref{eq-reuse-lemma}, we have
	\[
	\sup_{n\geq n_0}\frac{\sum _{k=0}^{Jn}|\widehat{B^n}(k)|w_{dk}}{w_n} \lesssim  \sup_{n\geq n_0} q^nw_n^{-1} <\infty,
	\]
	where the last inequality comes from $\liminf _n w_n ^{1/n}\geq 1$.
	
	Moreover, by \cite{Kahane1956} (see also \cite{BlyudzeShimorin1993}), $\Vert B^n\Vert _{\ell_a^1}\lesssim \sqrt{n}$, hence, by \eqref{assumption-weight-bddness} and the choice of $d$,
	\[
	\sup_{n\geq n_0}\frac{\sum _{k>Jn}|\widehat{B^n}(k)|w_{dk}}{w_n} \lesssim \sup_{n\geq n_0} \sqrt{n} \frac{w_{d\lfloor Jn\rfloor }}{w_n} \leq \sup_{n\geq n_0} \sqrt{n} \frac{w_{\mu n}}{w_n} < \infty.
	\]
	
	For the case $p=\infty$, let $f=\sum_{n\geq 0} a_n z^n \in \ell_a^{\infty}(w)$ and write
	\[
	C_{B_d}(f)=\sum_{n\geq 0} a_n \sum_{k\geq 0} \widehat{B^n}(k)z^{dk}=\sum_{k\geq 0} \left(\sum_{n\geq 0}a_n\widehat{B^n}(k)\right)z^{dk}.
	\]
	We have
	\[
	\left|\sum_{n\geq 0}a_n\widehat{B^n}(k)\right| \leq \sum_{n\geq 0}\left|a_nw_n\frac{\widehat{B^n}(k)}{w_n}\right| \leq \Vert f\Vert _{\ell_a^{\infty}(w)} \sum_{n\geq 0}\frac{\left|\widehat{B^n}(k)\right|}{w_n}.
	\]
	Then,
	\[
	\Vert f\circ B_d \Vert_{\ell_a^{\infty}(w)} \leq \Vert f\Vert _{\ell_a^{\infty}(w)} \sup_{k\geq0}\sum_{n\geq 0}\left|\widehat{B^n}(k)\right|\frac{w_{dk}}{w_n}.
	\]
	Now, by \eqref{eq-reuse-lemma}, we obtain
	\[
	\sum_{n>k/J}\left|\widehat{B^n}(k)\right|\frac{w_{dk}}{w_n} \leq C_0w_0\sum_{n>k/J} \frac{q^n}{w_n}<\infty,
	\]
	since $\liminf_n w_n^{1/n}\geq 1$. Now, set $m_k=\left\lfloor k/J\right\rfloor$. Since $w$ is decreasing, for $n\leq m_k$ we have $w_n\geq w_{m_k}$. Choosing $d$ large enough so that $dk\geq \mu m_k$ for every $k\geq1$, we obtain
	\[
	\sum_{n\leq k/J}\left|\widehat{B^n}(k)\right|\frac{w_{dk}}{w_n}\leq(m_k+1)\frac{w_{dk}}{w_{m_k}}\leq (m_k+1)\frac{w_{\mu m_k}}{w_{m_k}}\lesssim 1,
	\]
	where the last estimate follows from \eqref{assumption-weight-bddness}.
	
	\medskip
	
	The proof of \eqref{tiem2thmlast} goes along the same lines. Let $B$ be a finite Blaschke product vanishing at $0$. Writing $B(z)=z B_1(z)$ with $B_1$ a finite Blaschke product, we have, for every $d\in \N$,
	\[
	B^{dn}=z^{dn}B_1^{dn}=\sum_{{k\geq 0}}\widehat{B_1^{dn}}(k)z^{k+dn}.
	\]
	Thus, for every $d\geq 1$, one can check, as above, that $C_{B^d}$ is bounded on $\ell_a^1(w)$ and on $\ell_a^{\infty}(w)$ whenever, respectively,
	\begin{equation}\label{eqref-bounded}
		\sup_{n\geq 0}\sum _{k\geq 0}|\widehat{B_1^{dn}}(k)|\frac{w_{k+dn}}{w_n}  <\infty\quad \text{and}\quad \sup_{m\geq 0} \sum_{n= 0}^{\lfloor m/d \rfloor}\left|\widehat{B_1^{dn}}(m-dn)\right|\frac{w_{m}}{w_n}<\infty 
	\end{equation}
	Note that the function $B_1^d$ is a finite Blaschke product. Therefore, since $w_{k+dn}\leq w_{dn}$, the first part of \eqref{eqref-bounded} is a direct consequence of the assumption on $w$ and the estimate $\sum_{k\geq 0} |\widehat{B_1^{dn}}(k)|\lesssim \sqrt{dn+1}$.
	
	For the second part of \eqref{eqref-bounded}, use Remark \ref{remark-non-so-weight} to observe that, for every $d$ large enough, $\sup_{n\geq 0} n^2 w_{dn}/w_n<\infty$, hence $\sum_{n\geq 0} w_{dn}/w_n<\infty$. Moreover, since $w$ is decreasing, we have $w_m\leq w_{dn}$ for $m\geq dn$. It follows
	\[
	\sum_{n= 0}^{\lfloor m/d \rfloor}\left|\widehat{B_1^{dn}}(m-dn)\right|\frac{w_{m}}{w_n}\leq \sum_{n= 0}^{\lfloor m/d \rfloor}\frac{w_m}{w_n} \leq \sum_{n\geq 0} \frac{w_{dn}}{w_n}<\infty,
	\]
	upon choosing $d$ large enough.
\end{proof}

We indicate that a more thorough study of composition operators on general weighted Beurling--Sobolev spaces, including the Hilbert setting, is provided in the forthcoming paper \cite{CharpentierEspoullierZarouf2026-1}.

\bibliographystyle{amsplain}
\bibliography{refs}

\end{document}